\documentclass[reqno,11pt]{article}
\usepackage{mathrsfs}
\usepackage{amssymb}
\usepackage[usenames,dvipsnames]{xcolor}
\usepackage{amsthm}
\usepackage{amsmath}
\usepackage{amsfonts}
\usepackage{enumerate}
\usepackage{txfonts}

\usepackage{graphicx}

 \theoremstyle{definition}

 \numberwithin{equation}{section}

\newtheorem{theorem}{Theorem}[section]
\newtheorem{lemma}[theorem]{Lemma}
\newtheorem{corollary}[theorem]{Corollary}
\newtheorem{proposition}[theorem]{Proposition}

\theoremstyle{definition}
\newtheorem{definition}[theorem]{Definition}

\theoremstyle{remark}
\newtheorem{remark}[theorem]{Remark}

\begin{document}
\title{Remarks on asymptotic behaviors of strong solutions to a viscous Boussinesq system}
\author{Shangkun WENG\footnote{Pohang Mathematics Institute,
Pohang University of Science and Technology. Hyoja-Dong San 31, Nam-Gu
Pohang, Gyungbuk 790-784, Korea. {\it Email: skweng@postech.ac.kr, skwengmath@gmail.com}}}
\date{Pohang Mathematics Institute, POSTECH}
\maketitle

\def\be{\begin{eqnarray}}
\def\ee{\end{eqnarray}}
\def\p{\partial}
\def\no{\nonumber}
\def\e{\epsilon}
\def\de{\delta}
\def\De{\Delta}
\def\om{\omega}
\def\Om{\Omega}
\def\f{\frac}
\def\th{\theta}
\def\la{\lambda}
\def\b{\bigg}
\def\al{\alpha}
\def\La{\Lambda}
\def\ga{\gamma}
\def\Ga{\Gamma}
\def\ti{\tilde}
\def\Th{\Theta}
\def\si{\sigma}
\def\Si{\Sigma}
\def\bt{\begin{theorem}}
\def\et{\end{theorem}}
\def\bc{\begin{corollary}}
\def\ec{\end{corollary}}
\def\bl{\begin{lemma}}
\def\el{\end{lemma}}
\def\bp{\begin{proposition}}
\def\ep{\end{proposition}}
\def\br{\begin{remark}}
\def\er{\end{remark}}
\def\bd{\begin{definition}}
\def\ed{\end{definition}}
\def\bpf{\begin{proof}}
\def\epf{\end{proof}}
\begin{abstract}
In this paper, we first address the space-time decay properties for higher order derivatives of strong solutions to the Boussinesq system in the usual Sobolev space. The decay rates obtained here are optimal. The proof is based on a parabolic interpolation inequality, bootstrap argument and some weighted estimates. Secondly, we present a new solution integration formula for the Boussinesq system, which will be employed to establish the existence of strong solutions in scaling invariant function spaces. We further investigate the asymptotic profiles and decay properties of these strong solutions. Our results recover and extend the important results in Brandolese and Schonbek (Tran. A. M.S. Vol 364, No.10, 2012, 5057-5090).
\end{abstract}

\begin{center}
\begin{minipage}{5.5in}

\end{minipage}
\end{center}

\section{Introduction}

The Boussinesq system takes the following form
\be\label{bou}
\begin{cases}
\p_t u+ u\cdot\nabla u+ \nabla P= \Delta u+ \theta e_3,\\
\text{div u}=0,\\
\p_t \theta+ u\cdot\nabla \theta= \Delta \theta.
\end{cases}
\ee
where $u:\mathbb{R}^+\times \mathbb{R}^3\to \mathbb{R}^3$ is the velocity field, the scalars $P: \mathbb{R}^+\times \mathbb{R}^3\to \mathbb{R}$ and $\theta: \mathbb{R}^+\times \mathbb{R}^3\to \mathbb{R}$ denote the pressure and the temperature of the fluid, respectively. Here $e_3=(0,0,1)^T$, where $T$ is the transpose. We will consider the Cauchy problem to the Boussinseq system by prescribing the initial data
\be\no
u(0,x)=u_0(x),\quad \theta(0,x)=\theta_0(x),
\ee
where $u_0$ is divergent free.

The Boussinesq system has received significant attention in the mathematical fluid mechanics community due to its close connection to the 3-D incompressible flow. By putting $\theta\equiv 0$, we obtain the Navier-Stokes equations. The global existence of weak solutions to (\ref{bou}) do exist, however the uniqueness is still out of reach for current mathematical analysis. Strong solution with small initial data in different suitable function spaces are also been studies extensively. One may refer to \cite{ah07}, \cite{cd80}, \cite{dp08}, \cite{dp09}, \cite{gy96na}, \cite{kp} and the reference therein for more details. For the Boussinesq system with partial viscosity or fractional diffusion, one may refer to \cite{chae06} and \cite{jmwz} and the reference therein.

Our paper is inspired by the understanding of the important results in Brandolese and Schonbek \cite{bs12}, in which they investigated how the variations of the temperature will affect the asymptotic behavior of the velocity field. For the weak solution of (\ref{bou}) with initial data $(u_0,\theta_0)\in L^2$, they showed the growth of the energy $\|u(t)\|_{L^2}$ is generic, i.e. if $\int_{\mathbb{R}^3} \theta_0(y) dy\neq 0$, then there exists positive constants $c<C$
\be\no
c(1+t)^{1/4}\leq \|u(t)\|_{L^2} \leq C(1+t)^{1/4},\quad t\gg 1.
\ee
Under the assumption $\int_{\mathbb{R}^3} \theta_0(y) dy=0$ and some other smallness conditions, they also improved the estimate of $\|u(t)\|_{L^2}$:
\be\no
\|u(t)\|_{L^2} \leq C(1+t)^{-1/4}.
\ee
They established the exact pointwise asymptotic profiles of solutions in the parabolic region $|x|\gg \sqrt{t}$, which indicated the different behavior of the solution when $|x_3|\to \infty$ or when $\sqrt{x_1^2+x_2^2}\to \infty$. This also enabled them to rigorously prove the lower bound of energy.

Based on the $L^2$ decay estimates made in \cite{bs12}, we further investigate the space-time decay rates for higher order derivatives of strong solutions to (\ref{bou}). Using a parabolic interpolation inequality originally proved in \cite{k01} and a bootstrap argument, we first obtain the large time decay for higher order norms. For the far field behavior in space, we can directly obtain some weighted $L^2$ estimates for $\theta$, and then regard $\theta$ as forcing terms in the weighted estimates of the vorticity $\omega=\text{curl }u$, which will produce the corresponding estimates for the velocity field. Our proof uses and extends several ideas developed in Kukavica and Torres's series of papers in \cite{k01}-\cite{k09}.

In the second part of this paper, we will develop a new solution integration formula for the Boussinesq system, from which we can get some new information about how the variation of the temperature will affect the asymptotic behavior of the velocity field. We will employ this integration formula to construct some strong solutions in a scaling invariant space which is slight different from the one in \cite{bs12}. We will further investigate the asymptotic profile and decay properties of these strong solutions. We can recover and extend some of the results made in \cite{bs12}. We will illustrate the delicate difference between our and their results in the following sections. This paper will be organized as follows. Section \ref{section 2} will present the results about the space-time decay properties of higher order norms. Section \ref{section 3} concerns the existence and asymptotic behaviors of strong solutions in some scaling invariant function spaces.

\section{Space-time decay of Strong solutions in Sobolev spaces}\label{section 2}

The existence of weak solutions to the Boussinesq system have been proved in \cite{cd80}. In \cite{bs12}, the authors have proved the existence of weak solution with suitable decay in the long time regime, which was described in the following.

\bt\label{bou21}(Theorem 2.2 in \cite{bs12})
{\it
\begin{description}
  \item{(a)} Let $(u_0,\theta_0)\in L_{\sigma}^2\times L^2$. There exists a weak solution $(\theta, u)$ of the Boussinesq system (\ref{bou}), continuous from $\mathbb{R}^+$ to $L^2$ with the weak topology, with data $(u_0,\theta_0)$ such that, for any $T>0$,
  \be\no
  u\in L^2(0,T; V)\cap L^{\infty}(0,T; L_{\sigma}^2),\quad \theta\in L^2(0,T; H^1)\cap L^{\infty}(0,T; L^2).
  \ee
  Under the additional condition $\theta_0\in L^1$, then
      \be\label{bou211}
      \|\theta(t)\|_{L^2}^2 \leq C(t+1)^{-\f 32},\quad \|u(t)\|_{L^2}^2 \leq C (t+1)^{\f 12}.
      \ee
  Moreover, if $\theta_0\in L^1\cap L^p$ for some $1\leq p< \infty$, then
  \be\no
  \|\theta(t)\|_{L^p} \leq C(p) (t+1)^{-\f 32(1-\f 1p)}.
  \ee
  \item{(b)} (The $\int_{\mathbb{R}^3} \theta_0(y) dy=0$ case.) We additionally assume $\int_{\mathbb{R}^3} |y|\theta_0(y) dy<\infty$ and $\int_{\mathbb{R}^3} \theta_0(y) dy=0$. Then there exists an absolute constant $\e_0>0$ such that if
      \be\label{bou212}
      \|\theta_0\|_{L^1} <\e_0,
      \ee
      then the weak solution $(u,\theta)$ satisfies, for some constant $C>0$ and all $t\in \mathbb{R}^+$,
      \be\label{bou213}
      \|\theta(t)\|_{L^2}^2 \leq C(1+t)^{-\f 52},\quad \quad\|u(t)\|_{L^2}^2\rightarrow 0\quad \text{as $t\rightarrow \infty$}.
      \ee
      Moreover, under the additional condition $u_0\in L^{\f 32}\cap L_{\sigma}^2$, we have
      \be\label{bou214}
      \|u(t)\|_{L^2}^2 \leq C(1+t)^{-\f 12}.
      \ee
\end{description}
}
\et

Based on these $L^2$ decay estimates, we will investigate the space-time decay properties for high order derivatives of strong solutions of (\ref{bou}) in Sobolev space. To our purpose, we also present the following existence and uniqueness of mild solutions in the scaling invariant function space $L_{\sigma}^3(\mathbb{R}^3)\times L^1(\mathbb{R}^3)$, where a decay estimate of $\|u\|_{L^{\infty}(\mathbb{R}^3)}$ is available.

\bt\label{bou22}
{\it
Let $u_0\in L_{\sigma}^3(\mathbb{R}^3)$ and $\theta_0\in L^1(\mathbb{R}^3)$. There exists a number $\e_0>0$ sufficiently small so that if $\|u_0\|_{L^3(\mathbb{R}^3)}+\|\theta_0\|_{L^1(\mathbb{R}^3)}\leq \e_0$, then (\ref{bou}) admits a strong solution $(u,\theta)$ to (\ref{bou}) with initial data $(u_0,\theta_0)$, which satisfies for $t>0$
\be\label{bou221}
\|\nabla^k u(t)\|_{L^q(\mathbb{R}^3)} &\leq& C \e_0 t^{-\f k2-\f 32(\f 13-\f 1q)},\ \ k=0,1\ \text{and}\  3\leq q\leq \infty;\\\label{bou222}
\|\nabla^k\theta(t)\|_{L^q(\mathbb{R}^3)} &\leq& C\e_0 t^{-\f k2-1-\f32(\f 13-\f 1q)},\ \ k=0,1\ \text{and}\ 3\leq q\leq \infty.
\ee
}
\et
This is just an extension of the well-known results for the Navier-Stokes equations proved in Kato \cite{kato}, one can also refer to \cite{hs14} for the half-space case. Hence we have the following decay estimate for $\|u(t)\|_{L^{\infty}(\mathbb{R}^3)}$:
\be\no
\|u(t)\|_{L^{\infty}(\mathbb{R}^3)}\leq C t^{-\f 12}.
\ee

For simplicity, we assume that the initial data $(u_0,\theta_0)$ belong to the Schwartz class $\mathcal{S}$, so that for any $a\geq 0$ and $b\in \mathbb{N}_0=\{0,1,2,\cdots\}$,
\be\no
\||x|^a \nabla^b u_0\|_{L^2}<\infty,\quad \||x|^a \nabla^b B_0\|_{L^2}<\infty.
\ee
Assume the strong solution $(u,\theta)$ to (\ref{bou}) satisfies
\be\label{bwt11}
\|u(t)\|_{L^2}&=& O(t^{-\gamma}),\quad \|\theta(t)\|_{L^2}= O(t^{-\mu}),\\\label{bwt12}
\|u(t)\|_{L^{\infty}}&=&O(t^{-\f 12}),
\ee
where we assume that $\mu=\gamma+1$ and $\gamma\geq 0$. In the next section, we will derive a new solution integration formula (\ref{b4}), from which we can see this assumption $\mu=\gamma+1$ is quite reasonable.

By Theorem \ref{bou21} and \ref{bou22}, the strong solution $(u,\theta)$ with small initial data $(u_0,\theta_0)$ in the sense of Theorem \ref{bou21} and \ref{bou22} will satisfy (\ref{bwt11})-(\ref{bwt12}) with $\gamma=\f 14, \mu=\f 54$. Our main results in this section is stated as follows.
\bt\label{bou26}
{\it
Let $(u,\theta)$ be the strong solution to (\ref{bou}) with initial data $(u_0,\theta_0)$ belonging to the Schwartz class $\mathcal{S}$. Assume that (\ref{bwt11})-(\ref{bwt12}) are satisfied. Then we have the following weighted estimates for $u$ and $\theta$:
\be\label{bou261}
\||x|^a \nabla^b u(t,\cdot)\|_{L^p} &=&O(t^{-\gamma+\f{a}{2}-\f{b}{2}-\f{3}{4}(1-\f{2}{p})})
\ee
for any $b\in \mathbb{N}_0$ and $0\leq a <b+\f{5}{2}$ and $2\leq p\leq \infty$;
\be\label{bou262}
\||x|^a \nabla^b \theta(t,\cdot)\|_{L^p} =O(t^{-\mu+\f a2-\f b2-\f 34 (1-\f 2p)})
\ee
for all $b\in \mathbb{N}_0$ and $a\geq 0$ and $2\leq p\leq \infty$. Furthermore, for the vorticity $\omega(t,x)=\text{curl u}(t,x)$, we have
\be\label{bou263}
\||x|^a \nabla^b \omega(t,\cdot)\|_{L^p} =O(t^{-\gamma+\f a2-\f b2-\f 12-\f 34 (1-\f 2p)})
\ee
for all $b\in \mathbb{N}_0$ and $a\geq 0$ and $2\leq p\leq \infty$.
}
\et

\br\label{r11}
{\it
We find that the spatial decay property of the temperature is stronger than that of the velocity field in the sense that there is no restriction on the exponent of the weight. This is basically due to the pressure term in the velocity equations. Note that the spatial decay of the voricity field is also much stronger than the velocity field.
}
\er

\br\label{r12}
{\it
One can relax the conditions on the initial data, i.e. there are constants $r>0$ and $k\in \mathbb{N}_0$, such that for all $0\leq a\leq r$ and $0\leq b\leq k$
\be\no
\||x|^a \nabla^b u_0(\cdot)\|_{L^2}<\infty, \quad \||x|^a \nabla^b \theta_0(\cdot)\|_{L^2} <\infty.
\ee
Then the conclusions in Theorem \ref{bou26} also hold with some obvious modification.
}
\er

\subsection{Preliminary}

Before we start to prove Theorem \ref{bou26}, we collect several lemmas for our need. These lemmas are employed in the weighted estimates for higher order derivatives.
\bl\label{bwp1}
{\it
Let $\alpha_0>1,\alpha_1<1, \alpha_2<1$ and $\beta_1, \beta_2<1$. Assume that a continuously differentiable function $F:[1,\infty)\rightarrow [0,\infty)$ satisfies
\be\no
F'(t) &\leq& C_0 t^{-\alpha_0} F(t) + C_1 t^{-\alpha_1} F(t)^{\beta_1}+ C_2 t^{-\alpha_2} F(t)^{\beta_2}+C_3 t^{\gamma_2-1},\quad t\geq 1\\\no
F(1)  &\leq& K_0
\ee
where $C_0,C_1,C_2,C_3, K_0\geq 0$ and $\gamma_i=\f{1-\alpha_i}{1-\beta_i}>0$ for $i=1,2$. Assume that $\gamma_1\geq \gamma_2$, then there exists a constant $C^*$ depending on $\alpha_0,\alpha_1,\beta_1,\alpha_2, \beta_2,K_0,C_i, i=1,\cdots,4$, such that $F(t)\leq C^* t^{\gamma_1}$ for $t\geq 1$.
}
\el

\bl\label{bwp3}
{\it
Let $\tau_0>0$ and assume that $F: [\tau_0,\infty)\rightarrow [0,\infty)$ satisfies $\sup_{\tau_0\leq \tau\leq A} F(\tau)<\infty$ for all $A>\tau_0$. If there exist $C_0>0$ and $\gamma\in \mathbb{R}$ such that
\be\label{bwp31}
\sup_{\f t2\leq \tau\leq t} F(\tau)^2 \leq C_0 t^{-2\gamma} + C_0 t^{-\gamma} \sup_{\f t4\leq \tau\leq t} F(\tau),\quad t\geq 4\tau_0
\ee
then $F(t)=O(t^{-\gamma})$ as $t\rightarrow \infty$.
}
\el
Lemma \ref{bwp1} was proved in \cite{weng}, one can refer to \cite{kt07} for a simple version. Lemma \ref{bwp3} is just Lemma 3.2 in \cite{kt06}.

\bl\label{bwp2}
{\it
Let $p\in [1,\infty]$ and $T>0$. Assume that $u\in L^{\infty}(0,T: L^p(\mathbb{R}^n))$ and $t(u_t+(-\Delta)^{\alpha} u)\in L^{\infty}(0,T: L^p(\mathbb{R}^n))$. Then for $\alpha>\f 12$, there exists a constant $C=C(\alpha)>0$ such that
\be\no
\sup_{\f t2\leq \tau\leq t}\|\nabla u(\cdot,\tau)\|_{L^p} &\leq& C\b(\sup_{\f t4\leq\tau\leq t}\|u(\cdot,\tau)\|_{L^p}\b)^{\f{1}{2\alpha}}\b(\sup_{\f t4\leq\tau\leq t}\tau \|(u_t+(-\Delta)^{\alpha} u)(\cdot,\tau)\|_{L^p}\b)^{1-\f{1}{2\alpha}}\\\label{p4}
&\quad&+Ct^{-\f1{2\alpha}}\sup_{\f t4\leq\tau\leq t} \|u(\cdot,\tau)\|_{L^p}
\ee
holds for almost every $t\in (0,T)$.
}
\el

\bpf
This type of interpolation inequality has been proved in \cite{k01}. Here we generalize it to the fractional Laplacian case. The proof uses the solution integration formula. For $\forall t_0\in (0,T)$. We have
\be\no
u(t,x)= \int_{t_0}^t \int_{\mathbb{R}^3} G(t-s,x-y) f(y,s) dyds+ \int_{\mathbb{R}^3} G(t-t_0,x-y) u(y,t_0) dy
\ee
for a.e. $(t,x)\in (t_0,T)\times \mathbb{R}^n$, where $G(t,x)$ is the kernel of $(\p_t+(-\Delta)^{\alpha})u=0$. By Fourier transform, we get $\hat{G}(t,\xi)= e^{-t|\xi|^{2\alpha}}= e^{-|t^{\f{1}{2\alpha}}\xi|^{2\alpha}}$. Hence $G(t,x)= t^{-\f{n}{2\alpha}} K(x t^{-\f{1}{2\alpha}})$, where $K(x)$ satisfies $\hat{K}(\xi)=e^{-|\xi|^{2\alpha}}$. We can easily calculate the $L^p$ norm of $G(t,x)$ and $\nabla G(t,x)$:
\be\no
\|G(t,\cdot)\|_{L^p}= O(t^{-\f{n}{2\alpha}(1-\f 1p)}),\quad \|\nabla G(t,\cdot)\|_{L^p}= O(t^{-\f{1}{2\alpha}-\f{n}{2\alpha}(1-\f 1p)}),\quad 1\leq p\leq \infty.
\ee

By taking derivatives with respect to the space variables and estimating the $L^p$ norms, we get
\be\no
\|\nabla u(t)\|_{L^p} &\leq& C \int_{t_0}^t \|f(s)\|_{L^p}(t-s)^{-\f{1}{2\alpha}} ds + \|u(t_0)\|_{L^p}(t-t_0)^{-\f1{2\alpha}}\\\no
&\leq& C(t-t_0)^{1-\f{1}{2\alpha}} \sup_{t_0\leq s\leq t} \|f(s)\|_p + C\|u(t_0)\|_{L^p}(t-t_0)^{-\f{1}{2\alpha}},
\ee
from where
\be\label{p31}
\|\nabla u(t)\|_{L^p} &\leq& C(t-t_0)^{1-\f{1}{2\alpha}} \sup_{\f t2\leq s\leq t} \|f(s)\|_p + (t-t_0)^{-\f1{2\alpha}}\sup_{\f t2\leq s\leq t} \|u(s)\|_{L^p},
\ee
provided $0<\f t2\leq t_0\leq t$, where $C_0>0$ is a fixed constant depending only on $n$. Optimizing the right hand side of (\ref{p31}), we can conclude the result.

\epf

\subsection{Temporal decay rates for high order derivatives}

\bl\label{bwt}
{\it Assume that (\ref{bwt11}) and (\ref{bwt12}) hold with $\mu=\gamma+1$. Then we have the following temporal decay estimates for higher order derivatives: for $2\leq p\leq \infty$ and all multi-index $\alpha\in \mathbb{N}_0^3$
\be\label{bwt19}
\|\p_{\alpha} u(t)\|_{L^p}= O(t^{-\gamma-\f{|\alpha|}{2}-\f{3}{4}(1-\f 2p)}),\quad \|\p_{\alpha} \theta(t)\|_{L^p}= O(t^{-\mu-\f{|\alpha|}{2}-\f{3}{4}(1-\f 2p)}).
\ee
}
\el

\bpf
Our proof is based on parabolic interpolation lemma \ref{bwp2} and a bootstrap argument, which is different from \cite{sw} and is a modification of the method developed in \cite{k01}. Our decay rates are consistent with the scaling. Since $(u,\theta)$ satisfies (\ref{bou}), we can derive
\be\no
\nabla P=\nabla (-\Delta)^{-1} \text{div} (u\cdot\nabla u-\theta e_3).
\ee
Hence $\|\nabla P(t)\|_{L^p} \leq \|(u\cdot\nabla) u(t)\|_{L^p}+ \|\theta(t)\|_{L^p}$ for all $1<p<\infty$. Then by Lemma \ref{bwp2}, one gets
\be\no
\sup_{\f t2\leq \tau\leq t} \|\nabla u(\tau)\|_{L^2}^2 &\leq& \sup_{\f t4\leq\tau\leq t}\|u(\tau)\|_{L^2} \sup_{\f t4 \leq\tau\leq t}\b(\|u\cdot\nabla u(\tau)\|_{L^2}+\|\nabla P(\tau)\|_{L^2}+\|\theta(\tau)\|_{L^2}\b)\\\no
&\quad&\quad\quad+\f{C}{t}\sup_{\f t4\leq\tau\leq t}\|u(\tau)\|_{L^2}^2\\\no
&\leq&\sup_{\f t4\leq\tau\leq t}\|u(\tau)\|_{L^2} \sup_{\f t4\leq\tau\leq t}\b(\|u(\tau)\|_{L^{\infty}}\|\nabla u(\tau)\|_{L^2}+\|\theta(\tau)\|_{L^2}\b)\\\no
&\quad&\quad\quad+\f{C}{t}\sup_{\f t4\leq\tau\leq t}\|u(\tau)\|_{L^2}^2\\\no
&\leq& O(t^{-\gamma-\f12})\sup_{\f t4\leq\tau\leq t}\|\nabla u(\tau)\|_{L^2}+ O(t^{-\gamma-\mu})+ O(t^{-2\gamma-1}),
\ee
\be\no
\sup_{\f t2\leq \tau\leq t}\|\nabla\theta(\tau)\|_{L^2}^2 &\leq& \sup_{\f t4\leq \tau\leq t}\|\theta(\tau)\|_{L^2} \sup_{\f t4\leq \tau\leq t} \|u\cdot\nabla \theta(\tau)\|_{L^2}+\f{C}{t}\sup_{\f t4\leq\tau\leq t}\|\theta(\tau)\|_{L^2}^2\\\no
&\leq& O(t^{-\mu-\f 12}) \sup_{\f t4\leq\tau\leq t}\|\nabla \theta(\tau)\|_{L^2}+ O(t^{-2\mu-1}),
\ee
by Lemma \ref{bwp3}, we obtain
\be\label{bwt13}
\|\nabla u(t)\|_{L^2}&=& O(t^{-\gamma-\f 12}),\\\label{bwt14}
\|\nabla\theta(t)\|_{L^2}&=& O(t^{-\mu-\f 12}).
\ee
By Sobolev embedding and interpolation, (\ref{bwt11})-(\ref{bwt12}) and (\ref{bwt13})-(\ref{bwt14}) imply that for $2\leq p\leq 6$, we have
\be\label{bwt15}
\|u(t)\|_{L^p}= O(t^{-\gamma-\f 34(1-\f 2p)}),\quad \|\theta(t)\|_{L^p}= O(t^{-\mu-\f 34(1-\f 2p)}).
\ee
By similar arguments as above, we can easily obtain for $2\leq p\leq 6$
\be\label{bwt16}
\|\nabla u(t)\|_{L^p}= O(t^{-\gamma-\f 12-\f 34(1-\f 2p)}),\quad \|\nabla\theta(t)\|_{L^p}= O(t^{-\mu-\f 12-\f 34(1-\f 2p)}).
\ee
This will feedback with the estimate for all $L^p$ with $6\leq p\leq \infty$:
\be\label{bwt17}
\|u(t)\|_{L^p}= O(t^{-\gamma-\f 34(1-\f 2p)}),\quad \|\theta(t)\|_{L^p}= O(t^{-\mu-\f 34(1-\f 2p)})
\ee
and
\be\label{bwt18}
\|\nabla u(t)\|_{L^p}= O(t^{-\gamma-\f 12-\f 34(1-\f 2p)}),\quad \|\nabla\theta(t)\|_{L^p}= O(t^{-\mu-\f 12-\f 34(1-\f 2p)}).
\ee

We have proved (\ref{bwt19})for the cases $\alpha=0$ and $|\alpha|=1$. Assume that (\ref{bwt19}) holds for any $|\alpha|\leq k$, where $k\in \mathbb{N}$, we will show that (\ref{bwt19}) holds for all $|\alpha|=k+1$. Note that $\p_{\alpha} u$ and $\p_{\alpha} \theta$ satisfy the following equations
\be\no
(\p_t-\Delta)\p_{\alpha} u &=&-\nabla \p_{\alpha} P- \p_{\alpha}(u\cdot \nabla u)+ \p_{\alpha}\theta e_3\\\no
&=&-\nabla \p_{\alpha}P-\sum_{0\leq\beta\leq \alpha} \p_{\beta} u\cdot\nabla\p_{\alpha-\beta} u+ \p_{\alpha} \theta e_3,\\\no
(\p_t-\Delta)\p_{\alpha} \theta&=& -\p_{\alpha} (u\cdot\nabla \theta)=-\sum_{0\leq \beta\leq \alpha}C(\alpha,\beta) \p_{\beta} u\cdot\nabla \p_{\alpha-\beta} \theta.
\ee
Applying Lemma \ref{bwp2}, we obtain
\begin{equation}\no
\begin{aligned}
\sup_{\f t2\leq \tau\leq t}\|\nabla \p_{\alpha} u(\tau)\|_{L^p}^2 &\leq C\sup_{\f t4\leq\tau\leq t}\|\p_{\alpha} u(\tau)\|_{L^p}\sup_{\f t4\leq\tau\leq t}\b(\sum_{0\leq\beta\leq \alpha}\|\p_{\beta} u\cdot\nabla\p_{\alpha-\beta}u(\tau)\|_{L^p}^2\\\no
&\quad+\|\nabla \p_{\alpha} P(\tau)\|_{L^p}+\|\p_{\alpha}\theta(\tau)\|_{L^p}\b)+\f{C}{t}\sup_{\f t4\leq \tau\leq t} \|\p_{\alpha}u(\tau)\|_{L^p}\\\no
&\leq O(t^{-\gamma-\f{|\alpha|}{2}-\f 34 (1-\f 2p)})\sup_{\f t4\leq \tau\leq t}\b(\|u(\tau)\|_{L^{\infty}}\|\nabla \p_{\alpha}u(\tau)\|_{L^p}+\|\p_{\alpha} \theta(\tau)\|_{L^p}\\\no
&\quad+ \sum_{0<\beta\leq \alpha}\|\p_{\beta} u(\tau)\|_{L^{2p}}\|\nabla\p_{\alpha-\beta} u(\tau)\|_{L^{2p}}\b)+\f{C}{t} \sup_{\f t4\leq \tau\leq t}\|\p_{\alpha}u(\tau)\|_{L^p}^2\\\no
%&\leq O(t^{-\gamma-\f{|\alpha|}{2}-\f 34(1-\f 2p)})\b(O(t^{-\f 12})\sup_{\f t4\leq \tau\leq t}\|\nabla \p^{\alpha}u(\tau)\|_{L^p}\\\no
%&\quad+\sum_{0<\beta\leq \alpha} O(t^{-\gamma-\f{|\beta|}{2}-\f 34(1-\f 1p)}) O(t^{-\gamma-\f{|\alpha|-|\beta|+1}{2}-\f 34(1-\f 1p)})\\\no
%&\quad+O(t^{-\mu-\f{|\alpha|}{2}-\f 34(1-\f 2p)})\b)+ O(t^{-2\gamma-|\alpha|-\f 32(1-\f 2p)-1})\\\no
&\leq O(t^{-\gamma-\f{|\alpha|}{2}-\f{1}{2}-\f 34(1-\f 2p)}) \sup_{\f t4\leq \tau\leq t}\|\nabla \p_{\alpha}u(\tau)\|_{L^p}+ O(t^{-2\gamma-|\alpha|-\f 32(1-\f 2p)-1}),
\end{aligned}
\end{equation}

\begin{equation}\no
\begin{aligned}
\sup_{\f t2\leq \tau\leq t}\|\nabla \p_{\alpha}\theta(\tau)\|_{L^p}^2 &\leq C\sup_{\f t4\leq \tau\leq t}\|\p_{\alpha} \theta(\tau)\|_{L^p} \sup_{\f t4\leq \tau\leq t} \sum_{0\leq \beta\leq \alpha} \|\p_{\beta} u\cdot\nabla \p_{\alpha-\beta} \theta(\tau)\|_{L^p}\\\no
&\quad+ \f{C}{t} \sup_{\f t4\leq \tau\leq t} \|\p_{\alpha}\theta(\tau)\|_{L^p}^2\\\no
&\leq C\sup_{\f t4\leq \tau\leq t}\|\p_{\alpha} \theta(\tau)\|_{L^p} \sup_{\f t4\leq \tau\leq t} \b(\|u(\tau)\|_{L^{\infty}}\|\nabla \p_{\alpha} \theta(\tau)\|_{L^p}\\\no
&\quad+\sum_{0<\beta\leq \alpha} \|\p_{\beta} u\|_{L^{2p}}\|\nabla\p_{\alpha-\beta} \theta(\tau)\|_{L^{2p}}\b)+ \f{C}{t} \sup_{\f t4\leq \tau\leq t} \|\p_{\alpha}\theta(\tau)\|_{L^p}^2\\\no
%&\leq O(t^{-\mu-\f{|\alpha|}{2}-\f 34(1-\f 2p)})\b(O(t^{-\f12})\sup_{\f t4\leq \tau\leq t}\|\nabla\p^{\alpha}\theta(\tau)\|_{L^p}\\\no
%&\quad+\sum_{0<\beta\leq \alpha} O(t^{-\gamma-\f{|\beta|}{2}-\f 34(1-\f 1p)}) O(t^{-\mu-\f{|\alpha-\beta|+1}{2}-\f 34(1-\f 1p)})\b)+ O(t^{-2\mu-|\alpha|-\f 32(1-\f 2p)-1})\\\no
&\leq O(t^{-\mu-\f{|\alpha|+1}{2}-\f 34(1-\f 2p)})\sup_{\f t4\leq \tau\leq t}\|\nabla\p_{\alpha}\theta(\tau)\|_{L^p}+O(t^{-2\mu-|\alpha|-\f 32(1-\f 2p)-1}).
\end{aligned}
\end{equation}
These imply
\be\no
\|\nabla \p_{\alpha} u(t)\|_{L^p}= O(t^{-\gamma-\f{|\alpha|+1}{2}-\f{3}{4}(1-\f 2p)}),\quad \|\nabla \p_{\alpha} \theta(t)\|_{L^p}= O(t^{-\mu-\f{|\alpha|+1}{2}-\f{3}{4}(1-\f 2p)}).
\ee
\epf

\subsection{Weighted $L^2$ estimates}

By the results in previous section, we may assume that for $2\leq p\leq \infty$
\be\label{bw1}
\|\p_{\alpha} u(t,\cdot)\|_{L^p}= O(t^{-\gamma-\f{|\alpha|}{2}-\f 34(1-\f{2}{p})}), \quad \|\p_{\alpha}\theta(t,\cdot)\|_{L^p}= O(t^{-\mu-\f{|\alpha|}{2}-\f{3}{4}(1-\f {2}{p})})\quad \text{as $t\rightarrow \infty$}.
\ee

First we observe that the weighted estimate for $\theta$ can be obtained directly under the assumptions (\ref{bw1})
\bl\label{bw2}
{\it
Under the assumption (\ref{bw1}), we have
\be\label{bw201}
\||x|^a \theta(t,\cdot)\|_{L^2}= O(t^{-\mu+\f{a}{2}})\quad \text{as $t\rightarrow \infty$}
\ee
for all $a\geq 0$.
}
\el
\bpf

Define $G(t)=\int_{\mathbb{R}^3} |x|^{2a} |\theta(t,x)|^2 dx$, then
\be\no
&\quad&\f{d}{dt} G(t)+ 2 \int_{\mathbb{R}^3} |x|^{2a} |\nabla \theta(t,x)|^2 dx\\\no
&=& -4a \int_{\mathbb{R}^3} |x|^{2a-2} \theta(t,x) x\cdot\nabla\theta(t,x) dx + 2a\int_{\mathbb{R}^3} |x|^{2a-2} |\theta(t,x)|^2 x\cdot u(t,x) dx\\\no
&:=& I+II,
\ee
where
\be\no
|I| &\leq& 4a \int_{\mathbb{R}^3} |x|^{2a-1} |\theta(t,x)||\nabla \theta(t,x)| dx \\\no
&\leq& \int_{\mathbb{R}^3} |x|^{2a} |\nabla \theta(t,x)|^2 dx + C(a) \int_{\mathbb{R}^3} |x|^{2a-2} |\theta(t,x)|^2 dx\\\no
&\leq& \int_{\mathbb{R}^3} |x|^{2a} |\nabla \theta(t,x)|^2 dx + C(a) G(t)^{\f{a-1}{a}} \|\theta(t)\|_{L^2}^{\f 2a}
\ee
\be\no
|II| &\leq& C(a)\int_{\mathbb{R}^3} |\theta(t,x)|^2 |x|^{2a-1} |u(t,x)| dx \leq C(a)\|u(t)\|_{L^{\infty}}\int_{\mathbb{R}^3}|x|^{2a-1} |\theta(t,x)|^2 dx\\\no
&\leq& C(a)\|u(t)\|_{L^{\infty}} G(t)^{\f {2a-1}{2a}}\|\theta(t)\|_{L^2}^{\f 1a}.
\ee
Hence we have the following inequality:
\be\no
G'(t) &\leq& C(a)\|\theta(t)\|_{L^2}^{\f 2a} G(t)^{\f{a-1}{a}} + C(a)\|u(t)\|_{L^{\infty}} \|\theta(t)\|_{L^2}^{\f 1a} G(t)^{\f{2a-1}{2a}}\\\no
&\leq& C(a) t^{-\f{2\mu}{a}} G(t)^{\f{a-1}{a}} + C(a) t^{-\f 12-\f{\mu}{a}} G(t)^{\f{2a-1}{2a}}.
\ee
Now we apply Lemma \ref{bwp1}, where we take $C_0=C_3=0, C_1=C_2=C(a)$ and $\beta_1=\f{a-1}{a}<1, \alpha_1=\f{2\mu}{a}$, $\beta_2=\f{2a-1}{2a}<1, \alpha_2=\f 12+\f{\mu}{a}$. We first fix $a>2\mu$ so that $\alpha_1<1, \alpha_2<1$. Then $\gamma_1:=\f{1-\alpha_1}{1-\beta_1}=a-2\mu$ and $\gamma_2:=\f{1-\alpha_2}{1-\beta_2}=a-2\mu$. Hence Lemma \ref{bwp1} yields
\be\no
G(t)= O(t^{a-2\mu}), \quad \text{as $t\rightarrow \infty$}.
\ee
The conclusion for $a\in (0,2\mu]$ follows by interpolation.
\epf

Setting $\omega(t,x)=\text{curl u}(t,x)$, then $\omega$ satisfies
\be\no
\p_t \omega + u\cdot\nabla \omega-\omega\cdot\nabla u= \Delta \omega+ \text{curl} (\theta e_3).
\ee

\bl\label{bw3}
{\it
Under the assumptions (\ref{bw1}), we have the following estimate for all $a\geq 0$
\be\label{bw301}
\||x|^a\omega(t,\cdot)\|_{L^2}= O(t^{-\gamma-\f 12+ \f a2}).
\ee
}
\el

\bpf
Multiplying the vorticity equation by $2 |x|^{2a} \omega$ and setting $F(t)= \int_{\mathbb{R}^3} |x|^{2a} |\omega(t,x)|^2 dx$, then we get
\be\no
&\quad&\f{d}{dt} F(t) + 2\int_{\mathbb{R}^3} |x|^{2a} |\nabla \omega(t,x)|^2 dx \\\no
&=&- \int_{\mathbb{R}^3} 2|x|^{2a} \omega\cdot(u\cdot\nabla\omega) dx+ \int_{\mathbb{R}^3} 2|x|^{2a}\omega\cdot (\omega\cdot\nabla u) dx\\\no
&\quad&-4a \int_{\mathbb{R}^3} |x|^{2a-2} \sum_{i,j=1}^n x_j \omega_i \p_j\omega_i dx+ \int_{\mathbb{R}^3} 2|x|^{2a} \omega \cdot \text{curl}(\theta e_3) dx\\\no
&:=& I+ II + III+ IV.
\ee
These four terms will be estimated as follows.
\be\no
|I|&\leq& \f{1}{3}\int_{\mathbb{R}^3} |x|^{2a} |\nabla \omega|^2 dx+ C\|u\|_{L^{\infty}}^2 F(t), \quad |II|\leq 2\|\nabla u\|_{L^{\infty}} F(t),\\\no
|III|&\leq& C \int_{\mathbb{R}^3} |x|^{2a-1} |\omega| |\nabla \omega| dx \\\no
&\leq& \f{1}{3}\int_{\mathbb{R}^3} |x|^{2a} |\nabla \omega|^2 dx+ C F(t)^{\f{a-1}{a}} \|\omega\|_{L^2}^{\f{2}{a}},
\ee
\be\no
|IV|&=&2\b|\int_{\mathbb{R}^3} \text{curl}(|x|^{2a}\omega)\cdot (\theta e_3) dx\b|\\\no
&=& 2\b|\int_{\mathbb{R}^3} |x|^{2a} \text{curl} \omega \cdot (\theta e_3) dx+ \int_{\mathbb{R}^3}2a|x|^{2a-2}(x\times \omega)\cdot (\theta e_3) dx\b|\\\no
&\leq& C(a)\int_{\mathbb{R}^3} |x|^{2a} |\nabla \omega||\theta| dx + C(a) \int_{\mathbb{R}^3} |x|^{2a-1} |\omega| |\theta| dx\\\no
&\leq& \f{1}{3} \int_{\mathbb{R}^3} |x|^{2a} |\nabla\omega|^2 dx+ C(a)\int_{\mathbb{R}^3} |x|^{2a} |\omega|^2 dx+ C(a) F(t)^{\f{1}{2}} \||x|^{a-1}\theta\|_{L^2}.
\ee
Combining all these estimates together, we obtain
\be\no
F'(t)&\leq& C_0(\|u\|_{L^{\infty}}^2+\|\nabla u\|_{L^{\infty}})F(t)+ C_0 F(t)^{\f{a-1}{a}} \|\omega\|_{L^2}^{\f{2}{a}} \\\no
&\quad&\quad\quad + C_0 F(t)^{\f{1}{2}}  \||x|^{a-1} \theta\|_{L^2} \\\no
&\leq& C_0 t^{-\gamma-5/4}F(t)+ C_0t^{-\f{2}{a}(\gamma+\f{1}{2})} F(t)^{\f{a-1}{a}} +C_0 t^{-(\mu-\f{a-1}{2})} F(t)^{\f{1}{2}}.
\ee
Now we can apply Lemma \ref{bwp1}. Here $\beta_1=\f{a-1}{a}, \alpha_1=\f{2}{a}(\gamma+\f 12), \beta_2=\f{1}{2}, \alpha_2=\mu-\f{a-1}{2}$. To assure that $\alpha_1<1,\alpha_2<1$, we require $a>2(\gamma+\f 12)=(2\mu-1)$. Hence $\gamma_1=\f{1-\alpha_1}{1-\beta_1}=a-2(\gamma+\f 12)=\gamma_2=\f{1-\alpha_2}{1-\beta_2}=a-(2\mu+1)$. By Lemma \ref{bwp1}, we obtain
\be\no
F(t)\leq C t^{-2\gamma-1+a}.
\ee
\epf

By the relation $-\Delta u=\text{curl } \omega$ and the Caffarelli-Kohn-Nirenberg inequality \cite{ckn84}
\be\no
\||x|^a u\|_{L^p} \leq \||x|^{1+a} \nabla u\|_{L^p},
\ee
one can argue as in \cite{kt07} and \cite{k09} to obtain the weighted estimates for the velocity field $u$ as stated in the following theorem. Since the proof are almost the same, here we omit the details.
\bl\label{bw4}
{\it
Under the assumptions (\ref{bw1}), we have the following weighted estimates
\be\label{bw401}
\||x|^a u(t,\cdot)\|_{L^2}= O(t^{-\gamma_0+\f{a}{2}})
\ee
for all $a\in [0, \f{5}{2})$.
}
\el

\subsection{The Weighted estimates for higher order derivatives}

Based on the estimates (\ref{bw201}) and (\ref{bw401}), we can apply Lemma \ref{bwp2} to get the weighted estimates for higher order derivatives of $u$ and $\theta$.

\bl\label{bw11}
{\it
Under the assumptions  (\ref{bw1}), then the following estimates hold for all $a\geq 0, b\in \mathbb{N}_0$ and $2\leq p\leq\infty$
\be\label{bw111}
\||x|^a \nabla^b \theta(t,\cdot)\|_{L^p}= O(t^{-\mu-\f{b}{2}+\f{a}{2}-\f{3}{4}(1-\f{2}{p})}).
\ee
}
\el

\bpf
We only need to prove the case $a>2$, since we already know (\ref{bw111}) holds for $a=0$, the case $0<a\leq 2$ can be obtained by interpolation. By the Gagliardo-Nirenberg inequality, the case $p>2$ follows from $p=2$. Indeed, for any $f\in L^2(\mathbb{R}^3)\cap \dot{H}^2(\mathbb{R}^3)$, one has $\|f\|_{L^{\infty}(\mathbb{R}^3)}\leq \|f\|_{L^2(\mathbb{R}^3)}^{\f 14}\|f\|_{\dot{H}^2(\mathbb{R}^3)}^{\f 34}$. Hence one can derive the estimate of $\||x|^a \nabla^b \theta(t,\cdot)\|_{L^{\infty}(\mathbb{R}^3)}$ from those of $\||x|^{a} \nabla^b \theta(t,\cdot)\|_{L^2(\mathbb{R}^3)}$. The case $p\in (2,\infty)$ just follows from interpolation. Therefore, we assume $p=2$.
For $a>2$, we choose the weight $\phi$:
\be\no
\phi(t,x)= (|x|^2+t)^{\f{a}{2}},\quad t\geq 1,
\ee
then by simple calculations, we get
\be\no
|\nabla\phi(t,x)|\leq (|x|^2+t)^{\f{a-1}{2}},\quad |(\p_t-\Delta)\phi(t,x)|\leq C(|x|^2+t)^{\f{a}{2}-1}.
\ee

Assume that the conclusion holds for all the derivatives up to order $b\geq 0$, we want to show that it also holds for $b+1$. Take any $\alpha\in \mathbb{N}_0^3$ with $|\alpha|=b$, then
\be\no
(\p_t-\Delta)(\phi \p_{\alpha} \theta)&=& (\p_t-\Delta)\phi \p_{\alpha}\theta- 2\nabla\phi \cdot \nabla\p_{\alpha} \theta\\\no
&\quad&-\sum_{0\leq \beta\leq \alpha} C_{\alpha,\beta} \phi \p_{\beta} u\cdot \nabla\p_{\alpha-\beta} \theta.
\ee
Hence by Lemma \ref{bwp2} with $\alpha=1$, we obtain
\begin{equation}\no
\begin{aligned}
\sup_{\f t2\leq \tau\leq t}\|\nabla(\phi \p_{\alpha} \theta)\|_{L^2}^2 &\leq C\sup_{\f t4\leq \tau\leq t}\|\phi \p_{\alpha} \theta\|_{L^2} \sup_{\f t4\leq \tau\leq t}\b[\|(\p_t-\Delta)\phi \p_{\alpha}\theta\|_{L^2}+ \|\nabla \phi\cdot \nabla \p_{\alpha}\theta\|_{L^2}\\\no
&\quad+ \sum_{0\leq \beta\leq \alpha}C_{\alpha,\beta}\|\phi \p_{\beta}u\cdot \nabla\p_{\alpha-\beta}\theta\|_{L^2}\b]+\f{C}{t}\sup_{t/4\leq \tau\leq t} \|\phi \p_{\alpha}\theta(\tau)\|_{L^2}^2.
\end{aligned}
\end{equation}
By induction assumptions, we have
\be\no
\sup_{\f t2\leq \tau\leq t}\|\nabla(\phi \p_{\alpha} \theta)(\tau)\|_{L^2}^2&\geq&\f{1}{2}\sup_{\f t2\leq \tau\leq t}\|\phi\nabla\p_{\alpha} \theta\|_{L^2}^2-O(t^{-2\mu-b+a-1}),\\\no
\sup_{\f t4\leq \tau\leq t}\|\phi \p_{\alpha} \theta\|_{L^2}&\leq& O(t^{-\mu-\f{b}{2}+\f{a}{2}}),\quad\sup_{\f t4\leq \tau\leq t}\|(\p_t-\Delta)\phi \p_{\alpha} \theta\|_{L^2}\leq O(t^{-\mu-\f b2+\f a2-1}),\\\no
\sup_{\f t4\leq \tau\leq t}\|\nabla \phi\cdot \nabla\p_{\alpha} \theta\|_{L^2}&\leq& O(t^{-\f 12}) \sup_{\f t4\leq \tau\leq t}\|\phi\nabla \p_{\alpha} \theta\|_{L^2}.
\ee
For the last term, we estimate as follows.
\begin{equation}\no
\begin{aligned}
\|\phi\p_{\beta}u_j\p_j \p_{\alpha-\beta} \theta\|_{L^2} &\leq \begin{cases}
\|\p_{\beta} u_j\|_{L^{\infty}} \|\phi \p_j\p_{\alpha-\beta} \theta\|_{L^2},\quad &\text{if $|\beta|>0$}\\
\|u\|_{L^{\infty}} \|\phi \nabla \p_{\alpha} \theta\|_{L^2},\quad &\text{if $\beta=0$}
\end{cases}\\\no
&\leq\begin{cases}
O(t^{-\gamma-\mu-\f b2+\f a2-\f 12-\f 34}),\quad &\text{if $|\beta|>0$}\\
O(t^{-\gamma-\f 34})\|\phi \nabla \p_{\alpha} \theta\|_{L^2} ,\quad &\text{if $\beta=0$}.
\end{cases}
\end{aligned}
\end{equation}

Combining all these estimates together, we get
\be\no
\sup_{\f t2\leq\tau\leq t}\|\phi \nabla \p_{\alpha}\theta\|_{L^2}^2 &\leq& O(t^{-2\mu-b+a-1})+ O(t^{-\mu-\f b2+\f a2-\f 12})\sup_{\f t4\leq\tau\leq t}\|\phi \nabla \p_{\alpha}\theta\|_{L^2}.
\ee
This implies
\be\no
\|\phi \nabla \p_{\alpha} \theta(t,\cdot)\|_{L^2}= O(t^{-\mu-\f {b+1}2+\f a2}).
\ee

\epf

Now we show that the vorticity field has much stronger decay properties than the velocity field in the sense that there is no restriction on the exponent of the weight. This will help to improve the estimates on $u$ as shown in the following.
\bl\label{bw16}
{\it
Under the assumptions (\ref{bw1}), the following estimates
\be\label{bw161}
\||x|^a \nabla^b \omega(t,x)\|_{L^p}= O(t^{-\gamma-\f b2-\f 12+\f a2-\f{3}{4}(1-\f{2}{p})})
\ee
hold for any $a\geq 0$ and $b\in \mathbb{N}_0$ and $2\leq p\leq \infty$.
}
\el

\bpf
We choose same weight function as before. The conclusion is true for $b=0$ as showed in Lemma \ref{bw3}. We assume that the conclusion holds for any derivatives up to order $b$, we want to show that it also holds for $b+1$. Take any $\alpha\in \mathbb{N}_0^3$ with $|\alpha|=b$, then
\begin{equation}\no
\begin{aligned}
(\p_t-\Delta)(\phi \p_{\alpha}\omega)&=(\p_t-\Delta)\phi \p_{\alpha}\omega- 2(\nabla\phi\cdot\nabla)\p_{\alpha}\omega-\sum_{0\leq\beta\leq \alpha}C(\alpha,\beta)\b(\phi \p_{\beta} u\cdot\nabla\p_{\alpha-\beta}\omega\\\no
&\quad-\phi \p_{\beta}\omega\cdot\nabla\p_{\alpha-\beta}u\b)+ \phi \text{curl }(\p_{\alpha}\theta e_3).
\end{aligned}
\end{equation}
Then applying Lemma \ref{bwp2}, we obtain
\begin{equation}\no
\begin{aligned}
\sup_{\f t2\leq \tau\leq t}\|\nabla(\phi \p_{\alpha} \omega)\|_{L^2}^2 &\leq C\sup_{\f t4\leq \tau\leq t}\|\phi \p_{\alpha} \omega\|_{L^2} \sup_{\f t4\leq \tau\leq t}\b[\|(\p_t\phi-\Delta\phi) \p_{\alpha}\omega\|_{L^2}+ \|\nabla \phi\cdot \nabla \p_{\alpha}\omega\|_{L^2}\\\no
&\quad+ \sum_{0\leq \beta\leq \alpha}C_{\alpha,\beta}\b(\|\phi \p_{\beta}u\cdot\nabla\p_{\alpha-\beta}\omega\|_{L^2}+ \|\phi \p_{\beta}\omega\cdot\nabla \p_{\alpha-\beta}u\|_{L^2}\\\no
&\quad+\|\phi \text{curl }(\p_{\alpha} \theta e_3)\|_{L^2}\b)\b]+\f{C}{t}\sup_{\f t4\leq \tau\leq t} \|\phi \p_{\alpha}\omega\|_{L^2}^2
\end{aligned}
\end{equation}
As above, we have
\be\no
\sup_{\f t2\leq \tau\leq t}\|\nabla(\phi \p_{\alpha} \omega)\|_{L^2}^2&\geq&\f{1}{2}\sup_{\f t2\leq \tau\leq t}\|\phi\nabla\p_{\alpha} \omega\|_{L^2}^2-O(t^{-2\gamma-b-1+a-1}),\\\no
\sup_{\f t4\leq \tau\leq t}\|\phi \p_{\alpha} \omega\|_{L^2}&\leq& O(t^{-\gamma-\f{b}{2}-\f 12+\f{a}{2}}),\quad\sup_{\f t4\leq \tau\leq t}\|(\p_t-\Delta)\phi \p_{\alpha} \omega\|_{L^2}\leq O(t^{-\gamma-\f b2-\f 12+\f a2-1}),\\\no
\sup_{\f t4\leq \tau\leq t}\|\nabla \phi\nabla\p_{\alpha} \omega\|_{L^2}&\leq& O(t^{-\f 12}) \sup_{\f t4\leq \tau\leq t}\|\phi\nabla \p_{\alpha} \omega\|_{L^2}.
\ee
For the other three terms, we estimate as follows
\begin{equation}\no
\begin{aligned}
\|\phi\p_{\beta}u\cdot\nabla\p_{\alpha-\beta} \omega\|_{L^2} &\leq \begin{cases}\begin{array}{ll}
\|\p_{\beta} u\|_{L^{\infty}} \|\phi \nabla\p_{\alpha-\beta} \omega\|_{L^2},\quad &\text{if $|\beta|>0$}\\
\|u\|_{L^{\infty}} \|\phi \nabla \p_{\alpha} \omega\|_{L^2},\quad &\text{if $\beta=0$}
\end{array}
\end{cases}\\\no
&\leq\begin{cases}\begin{array}{ll}
O(t^{-2\gamma-\f b2+\f a2-1-\f 34}),\quad &\text{if $|\beta|>0$}\\
O(t^{-\gamma-\f 34})\|\phi \nabla \p_{\alpha} \omega\|_{L^2} ,\quad &\text{if $\beta=0$}
\end{array}
\end{cases}\\\no
\|\phi\p_{\beta} \omega_j \p_j\p_{\alpha-\beta}u\|_{L^2}&\leq \|\phi\p_{\beta} \omega_j\|_{L^2}\|\p_j\p_{\alpha-\beta}u\|_{L^{\infty}} \\\no
&\leq O(t^{-2\gamma-\f b2+\f a2-\f 34-1}),\\\no
\|\phi \text{curl }(\p_{\alpha}\theta e_3)\|_{L^2}&\leq O(t^{-\mu-\f b2-\f 12+\f a2}).
\end{aligned}
\end{equation}
Combining all these estimates together, we get
\be\no
\sup_{\f t2\leq\tau\leq t}\|\phi\nabla \p_{\alpha}\omega\|_{L^2}^2 &\leq& O(t^{-2\gamma-b+a-2})+ O(t^{-\gamma-\f b2+\f a2-1})\sup_{\f t4\leq\tau\leq t}\|\phi\nabla \p_{\alpha}\omega\|_{L^2}.
\ee
Hence $\|\phi \p_{\alpha}\omega\|_{L^2}= O(t^{-\gamma-\f{b+1}{2}-\f 12+\f a2})$.
\epf

In particular, we have showed that $\||x|^a u(t,\cdot)\|_{L^2}=O(t^{-\gamma+\f a2})$ for all $a\in [0,\f 52)$. Now one can argue as in Theorem 3.2 of \cite{kt06} to get the following theorem.
\bl\label{bw18}
{\it
Assume that (\ref{bw1}) hold, then
\be\label{bw181}
\||x|^a \nabla^b u(t,\cdot)\|_{L^p} = O(t^{-\gamma-\f b2+\f a2-\f 34(1-\f 2p)})
\ee
for every $2\leq p\leq \infty$, $a \in[0, b+\f{5}{2})$ and $b\in \mathbb{N}_0$.
}
\el

Finally, Theorem \ref{bou26} follows from Lemma \ref{bwt}, \ref{bw11}, \ref{bw16}, \ref{bw18}.

\section{Strong solutions in scaling invariant spaces}\label{section 3}

There is a natural scaling for (\ref{bou}), that is, for all $\lambda>0$, if $(u(t,x),\theta(t,x))$ is a smooth solution to (\ref{bou}) with initial data $(u_0(x), \theta_0(x))$, then
\be\no
u_{\la}(t,x)= \lambda u(\la^2 t, \la x),\quad \theta_{\la}(t,x) =\lambda^3 \theta(\la^2 t, \la x)
\ee
is also a smooth solution to (\ref{bou}) with initial data $(\lambda u_0(\la x), \la^3 \theta_0(\la x))$. We aim to establish the existence of strong solutions to (\ref{bou}) in suitable function space which is invariant under this scaling. Let us define the space $\mathcal{X}$ which consists of the vector function $u(t,x)$ equipped with the norm
\be\no
\|u\|_{\mathcal{X}}= \sup_{x\in \mathbb{R}^3, t>0}(\sqrt{t}+|x|) |u(t,x)|.
\ee
Denote the function space by $\mathcal{Y}$ which consists of the function $\theta(t,x)$ and equipped with the norm
\be\no
\|\theta\|_{\mathcal{Y}} = \sup_{x\in \mathbb{R}^3, t>0} (\sqrt{t}+|x|)^3 |\theta(t,x)|.
\ee

\bt\label{bou1}(Existence and uniqueness)
{\it
There exists an absolute constant $\e>0$ such that if the initial data $(u_0,\theta_0)$ satisfy the following conditions
\be\label{bou101}
\displaystyle\sup_{x\in \mathbb{R}^3} |x||u_0(x)|<\e,
\ee
and
\be\label{bou102}
\|t e^{t\Delta} \mathbb{P}(\theta_0 e_3)\|_{\mathcal{X}}<\e,\quad \ \|e^{t\Delta} \theta_0\|_{\mathcal{Y}}<\e,
\ee
where $u_0$ also satisfies the divergence-free condition. Then there exists a unique mild solution $(u,\theta)\in \mathcal{X}\times \mathcal{Y}$ to (\ref{bou}) such that
\be\no
\|u\|_{\mathcal{X}} \leq C\e\ \ \ \text{and}\ \ \|\theta\|_{\mathcal{Y}}\leq C\e,
\ee
where $C$ is a universal constant.
}
\et

\br\label{r1}
{\it
Here we remark that if $\theta_0\in L^1(\mathbb{R}^3)\cap \mathcal{Y}$ and $\|\theta_0\|_{L^1}+\|\theta_0\|_{\mathcal{Y}}<\e$, then as proved in \cite{bs12}, one can verify that
\be\no
\|t e^{t\Delta} \mathbb{P}(\theta_0 e_3)\|_{\mathcal{X}} + \|e^{t\Delta} \theta_0\|_{\mathcal{Y}}<\e.
\ee
Hence Theorem \ref{bou1} recovers the existence result presented in Proposition 2.4 in \cite{bs12}.
}
\er

\br\label{r2}
{\it Indeed, our result is broader than their results. Recall Theorem 1.2 (iii) in \cite{miya00}, where the following statement had been proved: Suppose there are functions $(a, b_1,b_2,b_3)(x)$ defined on $\mathbb{R}^3$, such that $a(x)=\sum_{k=1}^3 \p_k b_k(x)$, and
\be\no
|b_k(x)|\leq \e(1+|x|)^{-2},\ \ k=1,2,3, \ \ |a(x)|\leq C(1+|x|)^{-3},
\ee
then
\be\no
|(e^{t\Delta} a)(x)|\leq C \e(1+|x|)^{-\alpha} (1+t)^{-\beta/2}\quad \quad (\alpha\geq 0,\beta\geq 0, \alpha+\beta=3).
\ee
Back to our problem, we choose $b(x)= \f{\e x}{(1+|x|^2)^{\f 32}}$ and $\theta_0(x)= \sum_{k=1}^3 \p_k b_k= \f{3\e}{(1+|x|^2)^{3/2}}$, then we see that $\theta_0\not\in L^1(\mathbb{R}^3)$ but $e^{t\Delta} \theta_0\in \mathcal{Y}$ and satisfies (\ref{bou102}).
}
\er

Next we show that if the initial data have better decay rates, then we can improve the space time decay rates for the solution.
\bt\label{bou2}
{\it
\begin{description}
  \item{(i)} Let $u_0$ and $\theta_0$ be as in Theorem \ref{bou1}, and satisfy the additional decay estimates, for some $1\leq a<3,b\geq 3$, and a constant $C>0$,
      \be\label{bou201}
      \begin{cases}
      |u_0(x)|\leq C (1+|x|)^{-a},\\
      |\theta_0(x)|\leq C (1+|x|)^{-b}.
      \end{cases}
      \ee
      Then the solution constructed in Theorem \ref{bou1} satisfies, for another constant $C>0$ independent on $x$ and $t$,
      \be\label{bou202}
      |u(t,x)|&\leq& C\inf_{0\leq \eta\leq a} |x|^{-\eta} (1+t)^{(\eta-1)/2},\\\label{bou203}
      |\theta(t,x)|&\leq& C\inf_{0\leq \eta\leq b} |x|^{-\eta} (1+t)^{(\eta-3)/2}.
      \ee
  \item{(ii)} (The $\int_{\mathbb{R}^3} \theta_0(x) dx =0$ case.) Assume now $2\leq a<4, a\neq 3$ and $b\geq 4$, and let $u_0$ and $\theta_0$ satisfy the previous assumptions. If, in addition,
      \be\label{bou204}
      \int_{\mathbb{R}^3} \theta_0(x) dx=0,
      \ee
      then the decay of $u$ and $\theta$ is improved as follows:
      \be\label{bou205}
      |u(t,x)| &\leq& C\displaystyle\inf_{0\leq \eta\leq a}(1+|x|)^{-\eta}(t+1)^{(\eta-2)/2},\\\label{bou206}
      |\theta(t,x)| &\leq& C\displaystyle\inf_{0\leq \eta\leq b}(1+|x|)^{-\eta}(t+1)^{(\eta-4)/2}.
      \ee
\end{description}
}
\et

In the following, we first present a new integration formula for the Boussinesq system (\ref{bou}), and then prove Theorem \ref{bou1} and (\ref{bou2}).

\subsection{A new solution integration formula of (\ref{bou})}

In \cite{bs12} and \cite{kp}, they used the following integral formulation of the Boussinesq system
\be\label{formula1}
\begin{array}{llll}
u(t) &= e^{t\Delta} u_0- \int_0^t e^{(t-s)\Delta} \mathbb{P} \nabla \cdot(u\otimes u)(s) ds + \int_0^t e^{(t-s)\Delta} \mathbb{P}(\theta(s) e_3) ds\\
&:= e^{t\Delta} u_0 + B(u,u)+ L(\theta),\\
\theta(t)&= e^{t\Delta} \theta_0- \int_0^t e^{(t-s)\Delta} \nabla\cdot (\theta u)(s) ds\\
&:= e^{t\Delta} \theta_0 + \tilde{B}(\theta, u).
\end{array}
\ee

We present a new solution integration formula for the Boussinesq system.

\bl\label{integration formula}
{\it
Suppose $(u,\theta)$ is a smooth solution to the Boussinesq equation (\ref{bou}) with sufficient decay near the infinity, we have the following solution integration formula:
\be\label{integration}
\begin{cases}
\begin{array}{ll}
u(t,x) &=e^{t\Delta} u_0+ t e^{t\Delta} \mathbb{P}(\theta_0 e_3)-\int_0^t e^{(t-s)\Delta} \mathbb{P}\nabla\cdot(u\otimes u)(s) ds\\
&\quad -\int_0^t (t-s) e^{(t-s)\Delta} \mathbb{P}\nabla\cdot (\theta e_3\otimes u)(s)ds,\\
\theta(t,x) &= e^{t\Delta} \theta_0- \int_0^t e^{(t-s)\Delta} \nabla \cdot(\theta u)(s) ds\\
\end{array}
\end{cases}
\ee
}
\el

\bpf
We regard the nonlinear terms in (\ref{bou}) as a forcing term and rewrite (\ref{bou}) as
\be\label{b1}
\begin{cases}
\p_t u+\nabla P= \Delta u+ \theta e_3+ f(t,x),\\
\text{div u}=0,\\
\p_t \theta =\Delta \theta+ g(t,x),\\
u(0,x)=u_0(x),\quad \theta(0,x) =\theta_0(x),
\end{cases}
\ee
where $f(t,x)=-\nabla\cdot(u\otimes u)(t,x)$ and $g(t,x)= \nabla \cdot (\theta u)(t,x)$.

Taking the Fourier transform, we get
\be\label{b2}
\begin{cases}
\p_t \hat{u}(t,\xi)+ |\xi|^2 \hat{u}(t,\xi)+ i\xi \hat{P}(t,\xi)= \hat{\theta}(t,\xi) e_3+ \hat{f}(t,\xi),\\
\xi\cdot \hat{u}(t,\xi)=0,\\
\p_t\hat{\theta}(t,\xi)+ |\xi|^2\hat{\theta}(t,\xi)= \hat{g}(t,\xi),\\
\hat{u}(0,\xi)=\widehat{u_0}(\xi),\quad \hat{\theta}(0,\xi)=\widehat{\theta_0}(\xi).
\end{cases}
\ee
Simple calculations imply
\be\no
&&\hat{\theta}(t,\xi)= e^{-t|\xi|^2} \widehat{\theta_0}(\xi)+ \int_0^t e^{-(t-s)|\xi|^2} \hat{g}(s,\xi) ds,\\\no
&&\hat{P}(t,\xi)=-\f{i\xi}{|\xi|^2}\hat{\theta}(t,\xi)-\f{i\xi\cdot \hat{f}(t,\xi)}{|\xi|^2}.
\ee
Substituting the functions $(\theta, P)$ by the above formula, then we obtain
\be\no
\begin{cases}
\p_t \hat{u}(t,\xi)+ |\xi|^2 \hat{u}(t,\xi)=\b(\hat{f}(t,\xi)-\f{\xi\cdot \hat{f}(t,\xi)}{|\xi|^2}\xi\b)\\
\quad\quad\quad+\b(e^{-t|\xi|^2}\widehat{\theta_0}(\xi)+\int_0^t e^{-(t-s)|\xi|^2} \hat{g}(s,\xi) ds\b)\b(e_3-\f{e_3\cdot\xi}{|\xi|^2}\xi\b),\\
\hat{u}(0,\xi)=\widehat{u_0}(\xi).
\end{cases}
\ee
Then we get the representation formula for $u$:
\be\no
\hat{u}(t,\xi)&=& e^{-t|\xi|^2}\widehat{u_0}(\xi)+t e^{-t|\xi|^2}\widehat{\theta_0}(\xi)\b(e_3-\f{e_3\cdot\xi}{|\xi|^2}\xi\b)+\int_0^t e^{-(t-s)|\xi|^2}\b(\hat{f}(s,\xi)-\f{\xi\cdot \hat{f}(s,\xi)}{|\xi|^2}\xi\b) ds\\\no
&\quad&+ \int_0^t \int_0^s e^{-(t-\tau)|\xi|^2} \hat{g}(\tau,\xi)d\tau ds\b(e_3-\f{e_3\cdot\xi}{|\xi|^2}\xi\b)\\\no
&=&e^{-t|\xi|^2}\widehat{u_0}(\xi)+t e^{-t|\xi|^2}\widehat{\theta_0}(\xi)\b(e_3-\f{e_3\cdot\xi}{|\xi|^2}\xi\b)+\int_0^t e^{-(t-s)|\xi|^2}\b(\hat{f}(s,\xi)-\f{\xi\cdot \hat{f}(s,\xi)}{|\xi|^2}\xi\b) ds\\\label{b3}
&\quad&+ \int_0^t (t-s) e^{-(t-s)|\xi|^2} \hat{g}(s,\xi) ds\b(e_3-\f{e_3\cdot\xi}{|\xi|^2}\xi\b).
\ee

Hence
\be\label{b4}
\begin{cases}
\begin{array}{ll}
u(t,x) &= e^{t\Delta} u_0+ t e^{t\Delta} \mathbb{P}(\theta_0 e_3)-\int_0^t e^{(t-s)\Delta} \mathbb{P}\nabla\cdot(u\otimes u)(s) ds\\
&\quad -\int_0^t (t-s) e^{(t-s)\Delta} \mathbb{P}(\nabla\cdot (\theta u)(s) e_3) ds,\\
&:= e^{t\Delta} u_0 + t e^{t\Delta} \mathbb{P} (\theta_0 e_3)+ B(u,u) + E(u,\theta),\\
\theta(t,x) &= e^{t\Delta} \theta_0- \int_0^t e^{(t-s)\Delta} \nabla \cdot(\theta u)(s) ds\\
&:= e^{t\Delta} \theta_0+ \tilde{B}(u,\theta).
\end{array}
\end{cases}
\ee
\epf

\br\label{r3}
{\it
Indeed, one can obtain (\ref{integration}) from (\ref{formula1}) by a double iteration, i.e. substitute $\theta$ in the first equation of (\ref{formula1}) by the second equation.
}
\er

\subsection{Existence and decay properties of strong solutions}

We will use the following well-known fixed point theorem to prove Theorem \ref{bou1}.
\bl\label{fixed point}
{\it
Let $\mathcal{X}$ and $\mathcal{Y}$ be two Banach spaces, let $B: \mathcal{X}\times \mathcal{X}\to \mathcal{X}$, $E: \mathcal{X}\times \mathcal{Y}\to \mathcal{X}$ and $\tilde{B}: \mathcal{Y}\times \mathcal{X}\to \mathcal{Y}$ be three bilinear maps satisfying the estimates
\be\no
\|B(u,v)\|_{\mathcal{X}}\leq \alpha_1\|u\|_{\mathcal{X}} \|v\|_{\mathcal{X}},\quad \|E(u,\theta)\|_{\mathcal{X}}\leq \alpha_2 \|u\|_{\mathcal{X}}\|\theta\|_{\mathcal{Y}},\quad \|\tilde{B}(\theta,u)\|_{\mathcal{Y}}\leq \alpha_3 \|u\|_{\mathcal{X}} \|\theta\|_{\mathcal{Y}}
\ee
for some positive constants $\alpha_i, i=1,2,3$. Then there exists a positive constant $\e=\e(\alpha_1,\al_2,\al_3)$ such that for any $(U,\Theta)\in \mathcal{X}\times \mathcal{Y}$ satisfying
\be\no
\|U\|_{\mathcal{X}}+\|\Theta\|_{\mathcal{Y}} \leq \e,
\ee
the system
\be\no
u= U+ B(u,u)+ E(u,\theta),\quad \theta=\Theta + \tilde{B}(\theta,u)
\ee
has a unique solution $(u,\theta)\in \mathcal{X}\times \mathcal{Y}$ satisfying $\|u\|_{\mathcal{X}}+\|\theta\|_{\mathcal{Y}}\leq C\e$ for some universal constant $C$.
}
\el

To our purpose, we introduce the following function spaces. Let $a\geq 1$. We define $\mathcal{X}_{a}$ as the Banach space of divergence-free vector fields $u=u(t,x)$, defined and measurable on $\mathbb{R}^3\times \mathbb{R}^+$, such that, for some $C>0$,
\be\no
|u(t,x)|\leq C\inf_{0\leq \eta\leq a} |x|^{-\eta} (1+t)^{(\eta-1)/2},\quad \text{for all}\ t>0,x\in \mathbb{R}^3.
\ee
The norm $\|u\|_{\mathcal{X}_a}$ will be the infimum over all the above constant $C$.  Let $b\geq 3$, we define the space $\mathcal{Y}_b$ of function $\theta$ satisfying the estimates
\be\no
|\theta(t,x)|\leq C\inf_{0\leq \eta\leq b} |x|^{-\eta} (1+t)^{(\eta-3)/2},\quad \text{for all}\ t>0,x\in \mathbb{R}^3.
\ee
Define $\widetilde{\mathcal{X}}_a$ as the Banach space of divergence free vector fields $u=u(x,t)$ such that, for some $C>0$,
\be\no
|u(t,x)|\leq C\inf_{0\leq \eta\leq a} |x|^{-\eta} (1+t)^{(\eta-2)/2},\quad \text{for all}\ t>0,x\in \mathbb{R}^3.
\ee
For $b\geq 4$ we define the space $\widetilde{\mathcal{Y}}_b$ of functions $\theta$ satisfying the estimates
\be\no
|\theta(t,x)|\leq C\inf_{0\leq \eta\leq b} |x|^{-\eta} (1+t)^{(\eta-4)/2},\quad \text{for all}\ t>0,x\in \mathbb{R}^3.
\ee
Note that our function spaces $\mathcal{Y}_b$ and $\widetilde{\mathcal{Y}}_b$ are slightly different from those in \cite{bs12} by removing the integrability conditions on $\theta$: $\theta\in L_t^{\infty}(L^1)$ or $\theta\in L_t^{\infty}(L_1^1)$. We need the following $L^p$ estimates, which has been proved in \cite{bs12} by using Lorentz space technique.
\be\label{b601}
\|u(s)\|_{L^{q_1}} &\leq& C \|u\|_{\mathcal{X}}s^{-\f 12+\f{3}{2q_1}}, \quad 3<q_1\leq \infty,\\\label{b602}
\|\theta(s)\|_{L^{q_2}} &\leq& C\|\theta\|_{\mathcal{Y}} s^{-\f 32+\f{3}{2q_2}},\quad 1<q_2\leq \infty,\\\label{b603}
\|u(s)\|_{L^{q_1}} &\leq& C \|u\|_{\mathcal{X}_a}(1+s)^{-\f 12+\f{3}{2q_1}}, \quad 3/a<q_1\leq \infty,\\\label{b604}
\|\theta(s)\|_{L^{q_2}} &\leq& C\|\theta\|_{\mathcal{Y}_b} (1+s)^{-\f 32+\f{3}{2q_2}},\quad 3/b<q_2\leq \infty,\\\label{b605}
\|u(s)\|_{L^{q_1}} &\leq& C \|u\|_{\widetilde{\mathcal{X}}_a}(1+s)^{-1+\f{3}{2q_1}}, \quad 3/a<q_1\leq \infty,\\\label{b606}
\|\theta(s)\|_{L^{q_2}} &\leq& C\|\theta\|_{\widetilde{\mathcal{Y}}_b} (1+s)^{-2+\f{3}{2q_2}},\quad 3/b<q_2\leq \infty.
\ee
Indeed, one can prove them as follows. We take (\ref{b604}) as an example, others are similar. Assume that $\|\theta\|_{\mathcal{Y}_b}=1$, then
\be\no
\int_{\mathbb{R}^3} |\theta(s,x)|^{q_2} dx &\leq& \int_{|x|\leq R} (1+s)^{-\f 32 q_2} dx + \int_{|x|> R} |x|^{-bp_2} (1+s)^{\f{(b-3)p_2}{2}} dx \\\no
&\leq& c_1 (1+s)^{-\f 32 q_2} R^3 + c_2 R (1+s)^{\f{(b-3)q_2}{2}} R^{-bq_2+3},\quad \text{if}\ b q_2>3,\\\no
&\leq& c_3 (1+s)^{-\f 32 q_2+\f 32},
\ee
where the last inequality is obtained by optimizing the second line (taking $R=c (1+s)^{\f 12}$).

Before we prove Theorem \ref{bou1} and \ref{bou2}, we need to prepare some lemmas, which shows the boundedness of the operators $B(u,v), E(u,\theta)$ and $\tilde{B}(\theta, u)$ in different functions spaces. We first introduce some notations. Let $\mathbb{K}(t,x)$ be the kernel of $e^{t\Delta} \mathbb{P}$, let $F(t,x)$ be the kernel of $e^{t\Delta} \mathbb{P} \text{div} (\cdot)$. Both $\mathbb{K}(t,\cdot)$ and $F(t,\cdot)$ belong to $C^{\infty}(\mathbb{R}^3)$ and they satisfy the scaling properties
\be\no
\mathbb{K}(t,x)= t^{-3/2} \mathbb{K}(1,x/\sqrt{t}),\quad F(t,x)= t^{-2} F(1, x/\sqrt{t}).
\ee

By Proposition 1 in \cite{b09}, we have the following decomposition for $\mathbb{K}$ and $F$:
\be\label{f313}
\mathbb{K}(t,x)&=& \mathcal{R}(x)+ |x|^{-3} \Psi(x/\sqrt{t}),\\\label{f314}
F(t,x)&=& \mathcal{F}(x)+ |x|^{-4}\tilde{\Psi}(x/\sqrt{t}),
\ee
where $\mathcal{R}=(\mathcal{R}_{j,k})$ and $\mathcal{F}=(\mathcal{F}_{j;h,k})$ was defined as
\be\label{f311}
\mathcal{R}_{j,k}(x)= \p_{x_j x_k} E(x),\quad \mathcal{F}_{j;h,k}(x)= \p_{x_hx_j x_k} E(x)
\ee
where $E(x)=\f{1}{4\pi |x|}$; $\Psi$ and $\tilde{\Psi}$ are smooth outside the origin and such that, for all $\alpha\in \mathbb{N}^d$, and $x\neq 0$,
\be\no
|\p^{\alpha} \Psi(x)|+ |\p^{\alpha} \tilde{\Psi}(x)| \leq C e^{-c|x|^2}.
\ee
Here $C$ and $c$ are positive constants, depending on $|\alpha|$ but not on $x$.
\bl\label{b10}
{\it
\begin{description}
  \item{(i)} Let $1\leq a<3$. For some constant $C>0$, depending only on $a$, we have
\be\label{b101}
\|B(u,v)\|_{\mathcal{X}_a} &\leq& C\|u\|_{\mathcal{X}}\|v\|_{\mathcal{X}_a}\\\label{b102}
\|B(u,v)\|_{\mathcal{X}_{(2a)_*}} &\leq& C\|u\|_{\mathcal{X}_a}\|v\|_{\mathcal{X}_a},
\ee
where $(2a)_*= \min\{2a,4\}$. Moreover,
\be\label{b103}
\|B(u,v)\|_{\mathcal{X}}\leq C\|u\|_{\mathcal{X}}\|v\|_{\mathcal{X}}.
\ee
  \item{(ii)} Let $a\geq 1, b\geq 3$. For some constant $C>0$, depending only on $a,b$, we have
\be\label{b111}
\|\tilde{B}(\theta,u)\|_{\mathcal{Y}_{b}} &\leq& C\|u\|_{\mathcal{X}}\|\theta\|_{\mathcal{Y}_b},\\\label{b112}
\|\tilde{B}(\theta,u)\|_{\mathcal{Y}_{a+b}} &\leq& C\|u\|_{\mathcal{X}_a}\|\theta\|_{\mathcal{Y}_b}.
\ee
Moreover,
\be\label{b113}
\|\tilde{B}(\theta,u)\|_{\mathcal{Y}}\leq C\|u\|_{\mathcal{X}}\|\theta\|_{\mathcal{Y}}.
\ee
\item{(iii)} Let $2\leq a<4$ and $b\geq 4$. Then, for some constant $C>0$,
\be\label{b116}
\|B(u,v)\|_{\widetilde{\mathcal{X}}_a}&\leq& C\|u\|_{\mathcal{X}} \|v\|_{\widetilde{\mathcal{X}}_a},\quad \\\label{b117}
\|\tilde{B}(\theta, u)\|_{\widetilde{\mathcal{Y}}_b} &\leq& C\|u\|_{\mathcal{X}}\|\theta\|_{\widetilde{\mathcal{Y}}_b}.
\ee
\end{description}
}
\el

\bpf
Since $\mathcal{X}_a$ and $\widetilde{\mathcal{X}}_a$ are same as those in \cite{bs12}, the boundedness of $B(u,v)$ are exactly same as those in \cite{bs12}. Note that our definition of $\mathcal{Y}_b$ and $\widetilde{\mathcal{Y}}_b$ is slightly different from those in \cite{bs12}. However, one can check that in the proof of the boundedness of $B(u,v)$ and $\tilde{B}(\theta,u)$ in \cite{bs12}, they did not use the conditions $\theta \in L^{\infty}_t(L^1)$ and $\theta \in L^{\infty}(L^1_1)$. So their proof can also be applied to our cases.

\epf

\bl\label{b21}
{\it
There exists a constant $C>0$, such that
\be\label{b211}
\|E(u,\theta)\|_{\mathcal{X}} \leq C \|u\|_{\mathcal{X}} \|\theta\|_{\mathcal{Y}}.
\ee
}
\el

\bpf

We rewrite $E(u,\theta)$ as
\be\no
E(u,\theta) &=& \int_0^t (t-s) e^{(t-s)\Delta} \mathbb{P} \nabla \cdot [(\theta(s) e_3)\otimes u(s)] ds \\\no
&=& \int_0^t (t-s) \int F(t-s,x-y) (\theta(s,y) e_3\otimes u(s,y)) dy ds,
\ee
where $F(t,x)$ is the kernel of the operator $e^{t\Delta} \mathbb{P} \nabla$. Note that
\be\no
|F(t,x)|\leq C |x|^{-\eta} t^{-(4-\eta)/2},\quad \text{for all $0\leq \eta\leq 4$}
\ee
for some constant $C>0$ and $F(t,x)= t^{-2} F(1, x/\sqrt{t})$. This imply the following estimates:
\be\no
\|F(t)\|_p \leq C t^{-2+\f{3}{2p}},\quad 1\leq p\leq \infty.
\ee
Then by (\ref{b601}) and (\ref{b602}), we have
\be\no
\|E(u,\theta)(t)\|_{L_x^{\infty}} &\leq& C \int_0^t (t-s) \|F(t-s)\|_{L^{6}} \|u(s)\|_{L^6} \|\theta(s)\|_{L^{3/2}} ds \\\no
&\leq& C\int_0^t (t-s) (t-s)^{-2+\f 14} s^{-\f 12+\f 14} s^{-\f 32+1} ds =C\int_0^t (t-s)^{-\f 34} s^{-\f 34} ds \leq C t^{-\f 12}.
\ee

It remains to establish a pointwise estimate in the region $\{(x,t): |x|\geq 2\sqrt{t}\}$.  Decompose
\be\no
E(u,\theta) &=& \int_0^t (t-s)\b( \int_{|y|\leq |x|/2} +  \int_{|y|\geq |x|/2}\b) F(t-s,x-y) (\theta(s,y)e_3\otimes u(s,y)) dy ds\\\no
&:=& I_1'+I_2'.
\ee
We estimate $I_1'$ and $I_2'$ as follows.
\be\no
|I_1'| &\leq& C\int_0^t (t-s) \int_{|y|\leq |x|/2} |x-y|^{-1} (t-s)^{-\f32} (s^{\f 12}+|y|)^{-4} dy ds\\\no
&\leq& C|x|^{-1}\int_0^t (t-s)^{-\f 12} s^{-\f 12} ds\leq C|x|^{-1},\\\no
|I_2'| &\leq& C\int_0^t (t-s) \int_{|y|\geq |x|/2} |F(t-s,x-y)| (s^{\f 12}+|y|)^{-4} dy ds\\\no
&\leq& C \int_0^t (t-s) (s^{\f 12}+|x|)^{-4} \int_{|y|\geq |x|/2} |F(t-s,x-y)| dy ds,\\\no
&\leq& C \int_0^t (t-s)^{\f 12} (s^{\f 12}+|x|)^{-4} ds \leq C |x|^{-1} \int_{t/2}^{t} (t-s)^{\f 12} s^{-\f 32} ds + C |x|^{-4} \int_0^{t/2} (t-s)^{\f 12} ds\\\no
&\leq& C |x|^{-1} + C|x|^{-4} t^{\f 32}\leq C |x|^{-1},\quad \text{if $|x|\geq 2\sqrt{t}$}.
\ee

\epf

\bl\label{b26}
{\it
For $1\leq a<3$ and $b\geq 3$, there exists a constant $C$ depending only on $a$ and $b$, such that
\be\label{b261}
\|E(u,\theta)\|_{\mathcal{X}_a} &\leq& C \|u\|_{\mathcal{X}} \|\theta\|_{\mathcal{Y}_b},\\\label{b262}
\|E(u,\theta)\|_{\mathcal{X}_4} &\leq& C \|u\|_{\mathcal{X}_a} \|\theta\|_{\mathcal{Y}_b}.
\ee
}
\el

\bpf

Same as before, by (\ref{b601}) and (\ref{b604}), we have
\be\no
\|E(u,\theta)(t)\|_{L_x^{\infty}} &\leq& C \int_0^t (t-s) \|F(t-s)\|_{L^{6}} \|u(s)\|_{L^6} \|\theta(s)\|_{L^{3/2}} ds \\\no
&\leq& C\int_0^t (t-s) (t-s)^{-2+\f 14} s^{-\f 12+\f 14} (1+s)^{-\f 32+1} ds \\\label{bou265}
&=& C\int_0^t (t-s)^{-\f 34} s^{-\f 14} (1+s)^{-\f 12} ds \leq C (1+t)^{-\f 12}.
\ee

We estimate $I_1'$ and $I_2'$ as follows.
\be\no
|I_1'| &\leq& C\int_0^t (t-s) \int_{|y|\leq |x|/2} |x-y|^{-3} (t-s)^{-\f12} s^{-\f 12} |y|^{-a} (1+s)^{\f{a-3}{2}} dy ds\\\no
&\leq& C\int_0^t (t-s)^{\f 12} s^{-\f 12} (1+s)^{\f{a-3}{2}} ds |x|^{-3} \int_{|y|\leq |x|/2} |y|^{-a} dy,\\\no
&\leq& C |x|^{-a} (1+t)^{\f{a-1}{2}},\\\no
|I_2'| &\leq& C\int_0^t (t-s) \int_{|y|\geq |x|/2} |F(t-s,x-y)| s^{-\f 12} |y|^{-a} (1+s)^{\f{a-3}{2}} dy ds\\\no
&\leq& C |x|^{-a}\int_0^t (t-s) s^{-\f 12} (1+s)^{\f{a-3}{2}} \int_{|y|\geq |x|/2} |F(t-s,x-y)| dy ds,\\\no
&\leq& C |x|^{-a} \int_0^t (t-s)^{\f 12}s^{-\f 12} (1+s)^{\f{a-3}{2}} ds \leq C |x|^{-a} (1+t)^{\f{a-1}{2}}.
\ee
We have finished the proof of (\ref{b261}).

For the proof of (\ref{b262}), the estimate (\ref{bou265}) also holds. We estimate $I_1'$ and $I_2'$ as follows. By (\ref{b603}) and (\ref{b604}), for any fixed $p\in (\f 3a, \infty)$ and $q\in (1,\infty)$ satisfying $\f 1p+\f 1q=1$, we have
\be\no
|I_1'| &\leq& C\int_{0}^t (t-s) \int_{|y|\leq |x|/2} |x-y|^{-4} |u(s,y)||\theta(s,y)| dy ds\\\no
&\leq& C |x|^{-4} \int_0^t (t-s) \|u(s)\|_{L^p}\|\theta(s)\|_{L^q} ds\leq C |x|^{-4} \int_0^t (t-s) (1+s)^{-\f 12} ds\\\no
&\leq& C |x|^{-4} (1+t)^{\f 32},\\\no
|I_2'| &\leq& C \int_0^t (t-s) \int_{|y|>|x|/2} |K(t-s,x-y)||y|^{-1} |y|^{-3} dy\leq C|x|^{-4} \int_0^t (t-s) \|K(t-s)\|_{L^1} ds \\\no
&\leq& C |x|^{-4} (1+t)^{\f 32}.
\ee

\epf

\bl\label{b31}
{\it
For $2<a<4$ and $b\geq 4$, there exists a constant $C$ depending only on $a$ and $b$, such that
\be\label{b311}
\|E(u,\theta)\|_{\widetilde{\mathcal{X}}_a} \leq C \|u\|_{\mathcal{X}} \|\theta\|_{\widetilde{\mathcal{Y}}_b}.
\ee
}
\el

\bpf

By using (\ref{b601}) and (\ref{b606}), we have
\be\no
\|E(u,\theta)(t)\|_{L_x^{\infty}} &\leq& \int_0^{t/2} (t-s) \|F(t-s)\|_{L^{\infty}} \|u(s)\|_{L^4} \|\theta(s)\|_{L^{4/3}} ds+ \int_{t/2}^t (t-s)\|F(t-s)\|_{L^6} \|u(s)\|_{L^6} \|\theta(s)\|_{L^{3/2}} ds\\\no
&\leq& C\int_0^{t/2} (t-s)^{-1} s^{-\f 18} (1+s)^{-\f 78} ds+C\int_{t/2}^t (t-s)^{-\f 34} s^{-\f 14} (1+s)^{-1} ds \\\no
&\leq& C (1+t)^{-1}.
\ee

For $2<a<4$ and $b\geq 4$, we estimate $I_1'$ and $I_2'$ as follows.
\be\no
|I_1'| &\leq& C\int_0^t (t-s) \int_{|y|\leq |x|/2} |x-y|^{-\f{a+4}{2}} (t-s)^{-\f{4-\f{a+4}{2}}2} s^{-\f 12} |y|^{-\f{2+a}{2}} (1+s)^{\f{\f{2+a}{2}-4}{2}} dy ds\\\no
&\leq& C\int_0^t (t-s)^{{a}/4} s^{-\f 12} (1+s)^{a/4-3/2} ds |x|^{-\f{4+a}{2}} \int_{|y|\leq |x|/2} |y|^{-\f{2+a}{2}} dy,\\\no
&\leq& C |x|^{-a} (1+t)^{\f{a-2}{2}},
\ee
and
\be\no
|I_2'| &\leq& C\int_0^t (t-s) \int_{|y|\geq |x|/2} |F(t-s,x-y)| s^{-\f 12} |y|^{-a} (1+s)^{\f{a-4}{2}} dy ds\\\no
&\leq& C |x|^{-a}\int_0^t (t-s) s^{-\f 12} (1+s)^{\f{a-4}{2}} \int_{|y|\geq |x|/2} |F(t-s,x-y)| dy ds \quad \text{we need $a>2$},\\\no
&\leq& C |x|^{-a} \int_0^t (t-s)^{\f 12}s^{-\f 12} (1+s)^{\f{a-4}{2}} ds \leq C |x|^{-a} (1+t)^{\f{a-2}{2}}.
\ee
%\be\no
%|I_2'| &\leq& C\int_0^t (t-s) \int_{|y|\geq |x|/2} |F(t-s,x-y)| s^{-\f 12} |y|^{-\kappa} (1+s)^{\f{\kappa-4}{2}} dy ds\\\no
%&\leq& C\int_0^t (t-s)s^{-1/2}(1+s)^{\f{\kappa-4}{2}} \b(\int_{|y|\geq |x|/2} |F(t-s,x-y)|^p dy\b)^{\f 1p} \b(\int_{|y|\geq |x|/2} |y|^{-\kappa}dy\b)^{\f 1q} ds,\\\no
%&\leq& C |x|^{-\kappa+\f 3q} \int_0^t (t-s)\|F(t-s)\|_{L^p} s^{-\f 12} (1+s)^{\f{\kappa-4}{2}} ds \leq C |x|^{-a} (1+t)^{\f{a-1}{2}}.
%\ee
%Take $p=\f 32, q=3,\kappa=a+1$, then we have
%\be\no
%|I_2'| &\leq& C |x|^{-a} \int_0^t s^{-\f 12} (1+s)^{\f{a-3}{2}} ds \leq C |x|^{-a} (1+t)^{\f{a-2}{2}}.
%\ee
\epf

{\it Proof of Theorem \ref{bou1}.} The proof of Theorem \ref{bou1} will be based on the following iteration scheme: for $n=1,2,\cdots$, set
\be\label{bou105}
(u^0, \theta^0)&=&(U,\Theta)\coloneqq (e^{t\Delta} u_0+ t e^{t\Delta}\mathbb{P}(\theta_0 e_3), e^{t\Delta} \theta_0),\\\label{bou106}
(u^{n+1},\theta^{n+1})&=& (u^0+ B(u^n,u^n)+ E(u^n,\theta^n), \theta^0+ \tilde{B}(\theta^n, u^n)).
\ee
Note that $\|e^{t\Delta} u_0\|_{\mathcal{X}}\leq C \sup_{x\in \mathbb{R}^3} |x| |u_0(x)|$, then by (\ref{bou101}) and (\ref{bou102}), we have
\be\no
\|U\|_{\mathcal{X}} +\|\Theta\|_{\mathcal{Y}}\leq C\e.
\ee
By Lemma \ref{b10} and \ref{b21}, the conditions in Lemma \ref{fixed point} are satisfied, hence by Lemma \ref{fixed point}, there exists a unique mild solution $(u,\theta)$ to (\ref{bou}) with $\|u\|_{\mathcal{X}}+\|\theta\|_{\mathcal{Y}}\leq C\e$. We finished the proof.

{\it Proof of Theorem \ref{bou2}.} (i) As shown in \cite{b04,bs12}, with the assumption (\ref{bou201}), $e^{t\Delta} u_0, te^{t\Delta}\mathbb{P}(\theta_0 e_3)\in \mathcal{X}_a$ and $e^{t\Delta} \theta_0\in \mathcal{Y}_b$. Consider the approximate sequence $(u^n,\theta^n)$, by Lemma \ref{b10}, we have
\be\no
\|\theta^{n+1}\|_{\mathcal{Y}_b} \leq \|e^{t\Delta} \theta_0\|_{\mathcal{Y}_b} + C \e \|\theta^n\|_{\mathcal{Y}_b}.
\ee
If $\e>0$ is sufficiently small, so that $C\e<1$, then the sequence $(\theta^n)$ is bounded in $\mathcal{Y}_b$. Moreover, by Lemma \ref{b10} and \ref{b26}, we also have
\be\no
\|u^{n+1}\|_{\mathcal{X}_a} &\leq& \|e^{t\Delta} u_0\|_{\mathcal{X}_a}+ \|te^{t\Delta}\mathbb{P}(\theta_0 e_3)\|_{\mathcal{X}_a}+ C\e\|u^n\|_{\mathcal{X}_a}+ C\e \|\theta^n\|_{\mathcal{Y}_b},
\ee
hence $(u^{n+1})$ are also bounded in $\mathcal{X}_a$. Similarly, one can show $(u^n,\theta^n)$ is a Cauchy sequence in $\mathcal{X}_a\times \mathcal{Y}_b$. Thus the solution $(u,\theta)\in \mathcal{X}_a\times \mathcal{Y}_b$.

(ii) We need to verify that $(e^{t\Delta} u_0+ te^{t\Delta}\mathbb{P}(\theta_0 e_3))\in \widetilde{\mathcal{X}}_a$ and $e^{t\Delta} \theta_0\in \widetilde{\mathcal{Y}}_b$ under the assumptions in (ii). Here we only verify $te^{t\Delta}\mathbb{P}(\theta_0 e_3)\in \widetilde{\mathcal{X}}_a$, since other cases can be verified similarly. By (\ref{bou204}), we can rewrite
\be\no
t e^{t\Delta} \mathbb{P}(\theta_0 e_3) &=& t\int_{\mathbb{R}^3} \mathbb{K}(t,x-y)\theta_0(y) e_3 dy= t\int_{\mathbb{R}^3} (\mathbb{K}(t,x-y)-\mathbb{K}(t,x))\theta_0(y) e_3 dy\\\no
&=&t\int_{|y|\leq \sqrt{t}} (\mathbb{K}(t,x-y)-\mathbb{K}(t,x))\theta_0(y) e_3 dy + t\int_{|y|>\sqrt{t}} \mathbb{K}(t,x-y)\theta_0(y) e_3 dy\\\no
&\quad&\quad- t\mathbb{K}(t,x)\int_{|y|>\sqrt{t}} \theta_0(y) dy\\\no
&\eqqcolon & B_1+B_2 +B_3.
\ee
We can estimate $B_i, i=1,2,3$ as follows.
\be\no
|B_1|&=&\b|t\int_{|y|\leq \sqrt{t}} \int_0^1 (-y\cdot\nabla)\mathbb{K}(t,x-\lambda y) d\lambda \theta_0(y) e_3 dy\b|\\\no
&\leq&C tt^{-2} \int_{|y|\leq \sqrt{t}} |y|(1+|y|)^{-b} dy \leq C t^{-1},\\\no
|B_2| &\leq& t (1+\sqrt{t})^{-b} \int_{|y|>\sqrt{t}} |\mathbb{K}(t,x-y)|dy\leq  C (1+t)^{-b/2+1},\\\no
|B_3| &\leq& t t^{-3/2} \int_{|y|>\sqrt{t}} (1+|y|)^{-b} dy\leq C (1+t)^{-b/2+1}.
\ee
By the trivial estimate $|t e^{t\Delta}\mathbb{P}(\theta_0 e_3)|\leq C t$, we indeed obtain $|t e^{t\Delta}\mathbb{P}(\theta_0 e_3)|\leq C \min\{t, t^{-1}\}\leq C(1+t)^{-1}$. Similarly, we have
\be\no
t e^{t\Delta} \mathbb{P}(\theta_0 e_3)
&=&t\int_{|y|\leq |x|/2} (\mathbb{K}(t,x-y)-\mathbb{K}(t,x))\theta_0(y) e_3 dy + t\int_{|y|>|x|/2} \mathbb{K}(t,x-y)\theta_0(y) e_3 dy\\\no
&\quad&\quad- t\mathbb{K}(t,x)\int_{|y|>|x|/2} \theta_0(y) dy \eqqcolon A_1+A_2 +A_3.
\ee
We can estimate $A_i, i=1,2,3$ as follows.
\be\no
|A_2| &\leq& t |x|^{-b} \int_{|y|>|x|/2} |\mathbb{K}(t,x-y)|dy\leq t |x|^{-b}\|\mathbb{K}\|_{L^1}\leq C t|x|^{-b},\\\no
|A_3| &\leq& t |x|^{-3} \int_{|y|>|x|/2} (1+|y|)^{-b} dy\leq C t|x|^{-b},\\\no
|A_1|&=&\b|t\int_{|y|\leq |x|/2} \int_0^1 (-y\cdot\nabla)\mathbb{K}(t,x-\lambda y) d\lambda \theta_0(y) e_3 dy\b|\\\no
&\leq&C t|x|^{-4} \int_{|y|\leq |x|/2} |y|(1+|y|)^{-b} dy \leq C t |x|^{-4}.
\ee
In a word, $t e^{t\Delta} \mathbb{P}(\theta_0 e_3)\leq C t|x|^{-4}$, which implies $t e^{t\Delta}\mathbb{P}(\theta_0 e_3)\in \widetilde{\mathcal{X}}_a$ for $2\leq a\leq 4$. Then by Lemma \ref{b10} and \ref{b31}, one can prove that $(u^n,\theta^n)$ are uniformly bounded in $\widetilde{\mathcal{X}}_a\times \widetilde{\mathcal{Y}}_b$, which implies that $(u,\theta)\in \widetilde{\mathcal{X}}_a\times \widetilde{\mathcal{Y}}_b$.

\subsection{Asymptotic profiles of strong solutions}

\bt\label{bou5}
{\it
\begin{description}
  \item{(i)} Let $a>\f 32$ and $b\geq 3$. let $(u,\theta)$ be a mild solution of (\ref{bou}) satisfying the decay estimates (\ref{bou202})-(\ref{bou203}). Then the following profile for $u$ holds:
      \be\label{bou501}
      u(t,x) = e^{t\Delta} u_0(x) + t e^{t\Delta} \mathbb{P}(\theta_0 e_3) + \mathcal{R}_1(t,x),
      \ee
      where $\mathcal{R}_1(x,t)$ is a lower order term with respect to $t\nabla E_{x_3}(x)$ for $|x|\gg \sqrt{t}$, namely,
      \be\label{bou502}
      \lim_{\f{|x|}{\sqrt{t}}\rightarrow \infty} \f{\mathcal{R}_1(x,t)}{t |x|^{-3}}=0.
      \ee
      If $b>3$, we have for $j=1,2,3$:
      \be\label{bou5011}
      u_j(t,x) &=& e^{t\Delta} u_0(x) + t \mathbb{K}_{j,3}(t,x) \int_{\mathbb{R}^3} \theta_0(x) dx + \mathcal{R}_2(t,x)\\\no
      &=&e^{t\Delta} u_0(x)+t(\nabla E_{x_3})(x) \int_{\mathbb{R}^3} \theta_0(x)dx + \mathcal{R}_3(t,x),
      \ee
      where $\mathcal{R}_i(t,x), i=2,3 $ also share the same property (\ref{bou502}).
  \item{(ii)} (The $\int_{\mathbb{R}^3} \theta_0(x) dx=0$ case.) Assume now $a>2$ and $b>4$. Assume also that $\int_{\mathbb{R}^3} \theta_0(x) dx=0$. Let $(u,\theta)$ be a solution satisfying the decay condition (\ref{bou205}). Then the following profiles for $u_j (j=1,2,3)$ hold:
      \be\no
      u_j(x,t) &=& e^{t\Delta} u_0(x)- t\nabla \mathbb{K}_{j,3}(t,x)\cdot \int_{\mathbb{R}^3} y \theta_0(y) dy\\\no
      &\quad&\quad +\int_0^t (t-s) F(t-s,x)\cdot \int_{\mathbb{R}^3}(\theta e_3 \otimes u)(s,y) dy ds+\widetilde{\mathcal{R}}_1(t,x)\\\label{bou503}
      &=&e^{t\Delta} u_0(x)- t\nabla E_{x_j x_3}(x)\cdot\int_{\mathbb{R}^3} y \theta_0(y) dy\\\no
      &\quad&\quad+\nabla^2 \p_{x_j} E(x)\cdot\int_0^t(t-s)\int_{\mathbb{R}^3}(\theta e_3 \otimes u)(s,y) dy ds+ \widetilde{\mathcal{R}}_2(x,t),
      \ee
      where $\widetilde{\mathcal{R}}_i, i=1,2$ are lower order terms for $|x|\gg \sqrt{t}\gg 1$, namely
      \be\label{bou504}
      \lim_{t,\f{|x|}{\sqrt{t}}\rightarrow \infty} \f{\widetilde{\mathcal{R}}_i(x,t)}{t |x|^{-4}}=0,\quad \ i=1,2.
      \ee
\end{description}
}
\et

%\bc\label{bou6}
%{\it
%\begin{description}
%  \item{(a)} Let $a>\f 32, b>3$ and let $(u,\theta)$ be a solution as in part (a) of Theorem \ref{bou5}. Then for all $r,p$ such that
%  \be\no
%  r\geq 0,\quad 1<p<\infty, \quad r+\f 3p<\min\{a,3\},
%  \ee
%  there exists $t_0>0$ such that the solution satisfies the upper and lower estimates in the weighted $L^p$ norm
%  \be\label{bou601}
%  \phi(|m_0|) (1+t)^{\f 12(r+\f 3p-1)}\leq \|u(t)\|_{L_r^p}\leq C'(1+t)^{\f 12(r+\f 3p-1)}
%  \ee
%  for all $t\geq t_0$. Here $m_0=\int_{\mathbb{R}^3} \theta_0(x) dx$ and $\phi:\mathbb{R}^+\rightarrow \mathbb{R}$ is some continuous function such that $\phi(0)=0$ and $\phi(\sigma)>0$ if $\sigma>0$.
%  \item{(b)} ({\it The $\int_{\mathbb{R}^3} \theta_0(x) dx=0$ case.}) Under the assumptions of the previous item, with the stronger conditions $a>2, b>4$ and the additional zero mean condition $m_0=0$, let us set
%      \be\no
%      \tilde{m}=\liminf_{t\to\infty}\f{1}{t}\b|\int_{\mathbb{R}^3} y\theta_0(y)dy\b|
%      \ee
%       Then for all $r,p$ such that
%      \be\no
%      r\geq 0,\quad 1<p<\infty, \quad r+\f 3p<\min\{a,4\},
%      \ee
%      we have
%      \be\label{bou602}
%      \phi(|\tilde{m}|) (1+t)^{\f 12(r+\f 3p-2)} \leq \|u(t)\|_{L_r^p} \leq C'(1+t)^{\f 12(r+\f 3p-2)}
%      \ee
%      for another suitable continuous function $\phi:\mathbb{R}^+\rightarrow \mathbb{R}$ such that $\phi(0)=0$ and $\phi(\sigma)>0$ for $\sigma>0$.
%\end{description}
%}
%\ec

First we need the following asymptotic profile for $e^{t\Delta} \mathbb{P}(\theta_0 e_3)$.
\bl\label{bou8}
{\it
\begin{description}
  \item{(i)} Let $\theta_0(x)$ be any function satisfying the following estimate for some $3<b<4$
  \be\label{bou81}
  |\theta_0(x)|\leq C(1+|x|)^{-b}.
  \ee
  Then the $j$-component of $J(\theta_0)(t,x)\coloneqq e^{t\Delta} \mathbb{P}(\theta_0 e_3)$ can be decomposed as
  \be\no
  J(\theta_0)_j(t,x) &=& t \mathbb{K}_{j,3}(t,x) \int_{\mathbb{R}^3}\theta_0(x) dx + \mathcal{R}_j'(t,x)\\\label{bou801}
  &=& t (E_{x_j x_3})(x)\int_{\mathbb{R}^3}\theta_0(x) dx + \mathcal{R}_j^{''}(t,x),\quad \quad j=1,2,3,
  \ee
  where the remainder functions $\mathcal{R}'$ and $\mathcal{R}^{''}$ satisfy
  \be\label{bou802}
  |(\mathcal{R}',\mathcal{R}^{''})(t,x)| \leq C t |x|^{-b} ,\ \ \ \forall (t,x).
  \ee
  In particular, in the region $|x|\gg \sqrt{t}$, one has
  \be\no
  |(\mathcal{R}',\mathcal{R}^{''})(t,x)| \ll C t |E_{x_j,x_k}(x)|
  \ee
  along almost all directions.
  \item{(ii)} If $\theta_0(x)$ satisfies (\ref{bou81}) for some $4<b<5$ and $\int_{\mathbb{R}^3} \theta_0(y) dy=0$. Then the $j$-component $J(\theta_0)_j(t,x)$ can be decomposed as
      \be\no
      J(\theta_0)_j(t,x) &=& = -t \nabla \mathbb{K}_{j,3}(t,x) \cdot \b(\int_{\mathbb{R}^3} y \theta_0(y) dy\b) + \widetilde{\mathcal{R}}_j'(t,x)\\\label{bou805}
      &=& -t \nabla E_{x_jx_3}(x) \cdot \b(\int_{\mathbb{R}^3} y \theta_0(y) dy\b) + \widetilde{\mathcal{R}}_j^{''}(t,x),
      \ee
      where the remainder term $\widetilde{\mathcal{R}}$ satisfy
      \be\no
      |\widetilde{\mathcal{R}}(t,x)| &\leq& C t|x|^{-b},\quad \forall (t,x)\in \mathbb{R}^+\times \mathbb{R}^3.
      \ee
  \end{description}
}
\el

\bpf
(i). We first decompose $J(\theta_0)$ as follows:
\be\no
J(\theta_0)_j(t,x)&=& t \int_{\mathbb{R}^3} \mathbb{K}_{j,3}(x-y, t) \theta_0(y) dy\\\no
&=& t\b[\int_{|y|\leq |x|/2}+ \int_{|x-y|\leq |x|/2}+ \int_{|y|\geq |x|/2, |x-y|\geq |x|/2}\b] \mathbb{K}_{j,3}(x-y,t)\theta_0(y) dy\\\no
&\coloneqq& I_1+I_2+I_3.
\ee
$I_2$ and $I_3$ can be simply estimated as follows:
\be\no
|I_2(t,x)|&\leq& C t\int_{|x-y|\leq |x|/2} |\mathbb{K}_{j,3}(t,x-y)| (1+|y|)^{-b} dy \\\no
&\leq& C t|x|^{-b} \int_{\mathbb{R}^3} |\mathbb{K}(t,x-y)| dy\leq C t |x|^{-b},\\\no
|I_3(t,x)|&\leq& C t\int_{|x-y|\geq |x|/2, |y|\geq |x|/2} |x-y|^{-3} (1+|y|)^{-b} dy \\\no
&\leq& C t|x|^{-3} \int_{|y|\geq |x|/2} |y|^{-b} dy \leq C t |x|^{-b}.
\ee
It remains to treat $I_1$. We decompose $I_1$ as follows:
\be\no
I_1 &=& t\int_{\mathbb{R}^3} \mathbb{K}_{j,3}(t,x)\theta_0(y)dy- t\int_{|y|\geq |x|/2} \mathbb{K}_{j,3}(t,x)\theta_0(y) dy\\\no
&\quad&\quad+ t\int_{|y|\leq |x|/2}[\mathbb{K}_{j,3}(t,x-y)-\mathbb{K}_{j,3}(t,x)]\theta_0(y) dy\\\no
&\coloneqq& I_{1,1}+ I_{1,2} +I_{1,3},
\ee
where
\be\no
|I_{1,2}| &\leq& t |\mathbb{K}_{j,3}(t,x)|\int_{|y|\geq |x|/2} (1+|y|)^{-b} dy \leq C t |x|^{-b},\\\no
|I_{1,3}| &=& t \b|\int_{|y|\leq |x|/2} \int_0^1 \nabla\mathbb{K}_{j,3}(t,x-\lambda y)\cdot y d\lambda \theta_0(y) dy\b| \\\no
&\leq& t \int_{|y|\leq |x|/2}\max_{0\leq \la\leq 1} |\nabla \mathbb{K}(t,x-\lambda y)| |y| |\theta_0(y)| dy\\\no
&\leq& Ct|x|^{-4} \int_{|y|\leq |x|/2} |y| (1+|y|)^{-b} dy \leq C t|x|^{-b}.
\ee
Finally, for the estimate of $I_{1,1}$, we use the following decomposition of the kernel $\mathbb{K}$, established in \cite{b09}:
\be\no
\mathbb{K}_{j,k}(t,x)= E_{x_jx_k}(x) + |x|^{-3} \Psi_{j,k}(x/\sqrt{t}),\quad j,k=1,2,3,
\ee
where $\Psi_{j,k}$ is fast decaying: $|\Psi(y)|\leq C e^{-c|y|^2}$ for all $y\in \mathbb{R}^3$ and some constants $c, C>0$. Hence we can estimate $|\Psi(y)|\leq C |y|^{-b+3}$.

(ii). We argue as in (i). The estimates for $I_2, I_3$ and $I_{1,2}$ are the same, all of them are bounded by $C t |x|^{-b}$. However, we have $I_{1,1}=0$ by the vanishing condition. We need to further analyze $I_{1,3}$ to find the new main term. We can decompose $I_{1,3}$ as follows
\be\no
I_{1,3}&=& -t\int_{\mathbb{R}^3} \nabla \mathbb{K}_{j,3}(t,x) \cdot y\theta_0(y) dy + t\int_{|y|\geq |x|/2} \nabla \mathbb{K}_{j,3}(t,x) \cdot y\theta_0(y) dy\\\no
&\quad&\quad+t\int_{|y|\leq |x|/2}[\mathbb{K}_{j,3}(t,x-y)-\mathbb{K}_{j,3}(t,x)-(\nabla \mathbb{K}_{j,3}(t,x))\cdot(-y)]\theta_0(y) dy\\\no
&\coloneqq& I_{1,3,1}+ I_{1,3,2}+ I_{1,3,3}.
\ee
We can estimate $I_{1,3,2}$ and $I_{1,3,3}$ as follows
\be\no
|I_{1,3,2}| &\leq& C t |x|^{-4} \int_{|y|\geq |x|/2} |y|(1+|y|)^{-b} dy \leq C t |x|^{-b},\\\no
|I_{1,3,3}| &\leq& C t \int_{|y|\leq |x|/2} \max_{0\leq\la\leq 1} |\nabla_x^2 \mathbb{K}(t,x-\la y)| |y|^2 |\theta_0(y)| dy\\\no
&\leq& C t |x|^{-5}\int_{|y|\leq |x|/2} |y|^{2-b} dy \leq C t |x|^{-b}.
\ee
For $I_{1,3,1}$, we can argue as in (a) to get
\be\no
I_{1,3,1}(t,x) &=& -t \nabla \mathbb{K}_{j,3}(t,x) \cdot\int_{\mathbb{R}^3} y\theta_0(y) dy \\\no
&=& -t \nabla E_{x_j,x_3}(x)\cdot \int_{\mathbb{R}^3} y\theta_0(y) dy + O(t|x|^{-b}).
\ee
\epf
{\it Proof of Theorem \ref{bou5}}. (i). We can write
\be\no
u(t,x)= e^{t\Delta} u_0(x)+ t e^{t\Delta} \mathbb{P}(\theta_0 e_3)(x) +B(u,u)(t,x)+ E(u,\theta)(t,x).
\ee
Since $u\in\mathcal{X}_a$, $B(u,u)\in \mathcal{X}_{(2a)_*}$ and we have
\be\no
|B(u,u)(t,x)|\leq C |x|^{-(2a)_*} (1+t)^{((2a)_*-1)/2}, \quad \text{where}\ \ (2a)_*=\min\{2a,4\},
\ee
Since $a>\f 32$, we have $(2a)_*>3$ and hence
\be\no
\lim_{\f{|x|}{\sqrt{t}}\rightarrow \infty}\f{|B(u,u)(t,x)|}{t |x|^{-3}}\leq \lim_{\f{|x|}{\sqrt{t}}\rightarrow \infty} \b(\f{|x|}{\sqrt{t}}\b)^{-(2a)_*+3}=0.
\ee
Note that by Lemma \ref{b26}, we have $E(u,\theta)\in \mathcal{X}_4$, which also yields
\be\no
\lim_{\f{|x|}{\sqrt{t}}\rightarrow \infty}\f{|E(u,\theta)(t,x)|}{t |x|^{-3}}\leq \lim_{\f{|x|}{\sqrt{t}}\rightarrow \infty} \b(\f{|x|}{\sqrt{t}}\b)^{-(2a)_*+3}=0.
\ee

Moreover, Lemma \ref{bou8} and the condition $b>3$ guarantee that
\be\no
te^{t\Delta}\mathbb{P}(\theta_0 e_3)(x) =\b(\int_{\mathbb{R}^3} \theta_0(x) dx\b) t\nabla E_{x_3}(x)+ o(t|x|^{-3}),\quad \text{as $\f{|x|}{\sqrt{t}}\rightarrow\infty$}.
\ee

(ii). Note that in this case $u\in \widetilde{\mathcal{X}}_2$. As was proved in \cite{bs12}, we have $|B(u,u)(t,x)| \leq C |x|^{-4} t^{1/2}$. Now we turn to estimate $E(u,\theta)(t,x)$.
\be\no
E(u,\theta)(t,x) &=& \int_0^t (t-s) \b(\int_{|y|\leq |x|/2}+ \int_{|y|>|x|/2}\b) F(t-s, x-y) (\theta e_3\otimes u)(s,y) dy ds\\\no
&:=& I_1 +I_2.
\ee
$I_2$ can be easily bounded as follows.
\be\no
|I_2|&\leq& C\int_0^t (t-s) \int_{|y|>|x|/2}|F(t-s, x-y)| |y|^{-2} |y|^{-2} (s+1)^{-1} dy ds \\\no
&\leq& C |x|^{-4} \int_0^t (t-s)\|F(t-s,\cdot)\|_{L^1(\mathbb{R}^3)} (s+1)^{-1} ds = C|x|^{-4} \int_0^t (t-s)^{1/2} (s+1)^{-1} ds\\\no
&=& C|x|^{-4} \int_0^{t/2} t^{1/2} (s+1)^{-1} ds+ C |x|^{-4}\int_{t/2}^t (t-s)^{1/2} (t+1)^{-1} ds \\\no
&\leq& C |x|^{-4} t^{1/2} \log(1+t).
\ee
To estimate $I_1$, we rewrite $I_1$ as
\be\no
I_1(t,x) &=& \int_0^t (t-s) F(t-s,x)\int_{\mathbb{R}^3} (\theta e_3\otimes u)(s,y) dy ds\\\no
&\quad&\quad- \int_0^t (t-s) F(t-s,x) \int_{|y|>|x|/2} (\theta e_3\otimes u)(s,y) dy ds\\\no
&\quad& \quad+ \int_0^t (t-s) \int_{|y|\leq |x|/2} (F(t-s,x-y)-F(t,x)) (\theta e_3\otimes u)(s,y) dy ds\\\no
&:=& I_{1,1} + I_{1,2} +I_{1,3}.
\ee
Then $I_{1,2}$ and $I_{1,3}$ can be bounded in the following way
\be\no
|I_{1,2}|&\leq& C\int_0^t (t-s) |x|^{-3} (t-s)^{-1/2} \int_{|y|>|x|/2} |y|^{-2} |y|^{-2} (s+1)^{-1} dy ds \\\no
&\leq& C |x|^{-4} \int_0^t (t-s)^{1/2} (s+1)^{-1} ds \leq C |x|^{-4} t^{1/2} \log(1+t),\\\no
|I_{1,3}|&=&\b|-\int_0^t (t-s) \int_{|y|\leq |x|/2} \int_0^1 y\cdot\nabla F(t-s, x-\lambda y) d\lambda (\theta e_3\otimes u)(s,y) dy ds\b|\\\no
&\leq& C\int_0^t (t-s) \int_{|y|\leq |x|/2} \max_{0\leq \la \leq 1}|x-\la y|^{-4} (t-s)^{-1/2} |y| |u(s,y)| |\theta(s,y)| dy ds \\\no
&\leq& C|x|^{-4}\int_0^t (t-s)^{1/2} \int_{|y|\leq |x|/2} |y|(1+|y|)^{-2} (1+|y|)^{-5/2} (s+1)^{-3/4} dy ds\\\no
&\leq& C |x|^{-4} \int_0^t (t-s)^{1/2}(s+1)^{-3/4} ds \leq C|x|^{-4} t^{3/4}.
\ee

Applying Lemma \ref{bou8} will yield the asymptotic expansion (\ref{bou503}).

{\bf Acknowledgements.} The author would like to thank Prof. Jiahong Wu for reminding me of the interesting paper by Prof. Brandolese and Schonbek \cite{bs12}. Special thanks also go to Prof. Dongho Chae and Prof. Zhouping Xin for their encouragement and support.

\bibliographystyle{plain}

\begin{thebibliography}{10}

\bibitem{ah07}
H.Abidi, T. Hmidi. {\it On the global well posedness for Bounssinesq system}, J. Diff. Equ. 233, No.1, 199-220 (2007).

\bibitem{b04}
L. Brandolese. {\it Asymptotic behavior of the energy and pointwise estimates for solutions to the Navier-Stokes equations}. Rev. Mat. Iberoamericana 20 (2004), no. 1, 223--256.

\bibitem{b09}
L. Brandolese. {\it Fine properties of self-similar solutions of the Navier-Stokes equations.} Arch. Rational Mech. Anal. 192 (2009), 375-401.

\bibitem{bs12}
L. Brandolese, M. Schonbek. {\it Large time decay and growth for solutions of a viscous Boussinesq system.} Trans. Amer. Math. Soc. 364 (2012), no. 10, 5057--5090.

\bibitem{ckn84}
L. Caffarelli, R, Kohn, L. Nirenberg. {\it First order interpolation inequalities with weights.} Compositio Math. 53 (1984), 259-275.

\bibitem{cd80}
J.R. Cannon, E. DiBenedetto. {\it The initial problem for the Boussinesq equtions with data in $L^p$}, Approximation methods for Navier-Stokes problems (Proc. Sympos., Univ. Paderborn, Paderborn, 1979), pp. 129-144, Lecture Notes in Math., 771, Springer, 1980.

\bibitem{chae06}
D. Chae. {\it Global regularity for the 2D Boussinesq equations with partial viscosity terms.} Adv. in Math. 203, 497-513 (2006).

\bibitem{dp08}
R. Danchin, M. Paicu. {\it Existence and uniquness results for the Boussinesq system with data in Lorentz spaces.} Phys. D 237, N0.10-12, 1444-1460 (2008).

\bibitem{dp09}
R. Danchin, M. Paicu. {\it Global well-posedness issues for the inviscid Boussinesq system with Yudovich's type data.} Comm. Math. Phys. 290 (2009), no. 1, 1--14.

\bibitem{gy96na}
B. Guo, G. Yuan. {\it On the suitable weak solutions for the Cauchy problem of the Boussinesq equations.} Nonlinear Anal. 26 (1996), no. 8, 1367--1385.

\bibitem{hs14}
P. Han, M. Schonbek. {\it Large time decay properties of solutions to a viscous Boussinesq system in a half space.} Adv. Differential Equations 19 (2014), no. 1-2, 87--132.

\bibitem{jmwz}
Q. Jiu, C. Miao, J. Wu, and Z. Zhang. {\it The 2D incompressible Boussinesq equations with general
critical dissipation.} SIAM Journal on Mathematical Analysis, 46 (2014), No.5, 3426--3454.

\bibitem{kp}
G. Karch, N.Prioux. {\it Self-similarity in viscous Bounssinesq equation.} Proc. Amer. Math. Soc. 136, 879-888 (2008).

\bibitem{kato}
T. Kato. {\it Strong $L^p$-solutions of the Navier-Stokes equation in $\mathbb{R}^m$, with applications to weak solutions.} Math. Z. 187 (1984), 471-480.

\bibitem{k01}
I. Kukavica. {\it Space-time Decay for Solutions of the Navier-Stokes Equations.} Indiana University Mathematics Journal, vol. 50, no. 1(2001), 205-222.

\bibitem{kt06}
I. Kukavica; J. J. Torres. {\it Weighted bounds for the velocity and the vorticity for the Navier-Stokes equations.} Nonlinearity 19 (2006), no. 2, 293--303.

\bibitem{kt07}
I. Kukavica; J. J. Torres. {\it Weighted $L^p$ decay for solutions of the Navier-Stokes equations.} Comm. Partial Differential Equations 32 (2007), no. 4-6, 819--831.

\bibitem{k09}
I. Kukavica. {\it On the weighted decay for solutions of the Navier-Stokes system.} Nonlinear Anal. 70 (2009), no. 6, 2466--2470.

\bibitem{miya00}
T. Miyakawa.{\it On space-time decay properties of nonstationary incompressible Navier-Stokes flows in $\mathbb{R}^n$.} Funkcialaj Ekvacioj, 43 (2000), 541-557.

%\bibitem{st}
%O. Sawada, Y. Taniuchi. {\it On the Boussinesq flow with nondecaying initial data.} Funk. Ekvac. 47, No.2, 225-250 (2004).

\bibitem{sw}
M. Schonbek, M. Wiegner.{\it On the decay of higher-order norms of the solutions of Navier-Stokes equations.} Proc. Roy. Soc. Edinburgh Sect. A 126 (1996), no. 3, 677-685.

\bibitem{weng}
S. Weng. {\it Space-time decay estimates for the incompressible viscous resistive MHD and Hall-MHD equations.} arXiv:1405.4922.

\end{thebibliography}

\end{document}